\documentclass[opre,nonblindrev]{informs3}

\DoubleSpacedXI 

\usepackage{endnotes}
\let\footnote=\endnote

\usepackage{tikz}
\usepackage{bm}
\usepackage{comment}
\usepackage{subcaption}
\usepackage{multirow}							

\usepackage[ruled,vlined,linesnumbered]{algorithm2e}

\usepackage{mathtools}

\newcommand{\suchthat}{\;\ifnum\currentgrouptype=16 \middle\fi|\;}
\newcommand{\RU}{RU}
\newcommand{\SU}{SU}
\newcommand{\RD}{RD}
\newcommand{\SD}{SD}

\def\Z{\mathbb{Z}}
\def\Re{{\mathbb R}}

\def\conv{\mathop{\rm conv}}
\def\dim{\mathop{\rm dim}}

\def\proj{\mathop{\rm proj}}

\def\Re{{\mathbb R}}

\def\mc{\mathcal}	

\newcommand{\vc}[1]{\bm{#1}}	
\newcommand{\mr}[1]{\mathrm{#1}}	
\newcommand{\mt}[1]{\mathtt{#1}}	

\def\Sol{\mr{Sol}}

\newcommand{\blue}[1]{\textcolor{black}{#1}}

\usepackage{natbib}
 \bibpunct[, ]{(}{)}{,}{a}{}{,}%
 \def\bibfont{\small}%
\TheoremsNumberedThrough     
\ECRepeatTheorems

\EquationsNumberedThrough    

\begin{document}


\RUNAUTHOR{Salemi and Davarnia}

\RUNTITLE{On the Structure of DD-Representable MIPs with Application to UCP}

\TITLE{On the Structure of Decision Diagram-Representable Mixed Integer Programs with Application to Unit Commitment}

\ARTICLEAUTHORS{%
\AUTHOR{Hosseinali Salemi}
\AFF{Department of Industrial and Manufacturing Systems Engineering, Iowa State University, Ames, IA 50011, \EMAIL{hsalemi@iastate.edu}}
\AUTHOR{Danial Davarnia}
\AFF{Department of Industrial and Manufacturing Systems Engineering, Iowa State University, Ames, IA 50011, \EMAIL{davarnia@iastate.edu}} 
} 

\ABSTRACT{

Over the past decade, decision diagrams (DDs) have been used to model and solve integer programming and combinatorial optimization problems. 
Despite successful performance of DDs in solving various discrete optimization problems, their extension to model mixed integer programs (MIPs) such as those appearing in energy applications has been lacking.
More broadly, the question which problem structures admit a DD representation is still open in the DDs community. 
In this paper, we address this question by introducing a geometric decomposition framework based on rectangular formations that provides both necessary and sufficient conditions for a general MIP to be representable by DDs. 
As a special case, we show that any bounded mixed integer linear program admits a DD representation through a specialized Benders decomposition technique. 
The resulting DD encodes both integer and continuous variables, and therefore is amenable to the addition of feasibility and optimality cuts through refinement procedures. 
As an application for this framework, we develop a novel solution methodology for the unit commitment problem (UCP) in the wholesale electricity market.
Computational experiments conducted on a stochastic variant of the UCP show a significant improvement of the solution time for the proposed method when compared to the outcome of modern solvers.  
}


\KEYWORDS{Decision Diagrams, Benders Decomposition, Mixed Integer Programs, Unit Commitment} 

\maketitle

%


\section{Introduction} \label{sec:introduction}

Improved methods for scheduling and dispatching resources are necessary to accommodate the increases in renewable generation, distributed energy resources and storage in the electric grid.
At the core of electric grid lies the Unit Commitment Problem (UCP), which determines
the thermal unit schedules and generation dispatch in wholesale electricity markets.
Unfortunately, the UCP is computationally difficult, making it challenging for available software despite all advances in computing power and technology. 
For this reason, innovative methodologies that produce optimal solutions to the UCP satisfying realistic yet complicated constraints are needed to manage the increasing amounts of variable renewable generation and distributed energy resources in the electric grid.

Over the past decade, a powerful alternative to traditional solution techniques was developed based on a new paradigm called Decision Diagrams (DDs), in which discrete optimization problems are reconstructed as network models with special structure. 
The network structure is exploited to speed up the solution process and improve the solution quality, especially for formulations that are poorly handled by conventional solvers.
However, since DDs have been designed for modeling discrete problems, they have never been used to solve optimization problems that appear in the electric grid applications, as they naturally contain both discrete and continuous variables. 
In this paper, we generalize the concept of DDs to mixed integer programs (MIPs), and thereby, establish a novel framework that can be applied to the UCP with the purpose of improving the solution time and quality obtained from existing techniques.

\subsection{Background on Decision Diagrams} \label{subsec:background_DD}

DDs were introduced in 2006 \citep{hadzic:ho:2006} as a solution method for discrete optimization and combinatorial problems. 
DDs' principal idea is to model the solutions of the underlying discrete set through a directed acyclic graphical structure with a root and a terminal node.
Each path of the graph from the root to the terminal node encodes a solution of the problem where arcs are labeled by values of decision variables on their relative domain.
Structurally, DDs resemble branch-and-bound trees and dynamic programming graphs.
The key difference, however, is that the size of DDs can be controlled by a width limit factor that mitigates the exponential growth rate intrinsic to branch-and-bound and dynamic programming.
This control factor is attributed to a node-merging concept where nodes of the DD are merged to reduce its size at the price of expanding the solution set, thereby yielding a discrete relaxation.
By exploiting the concept of relaxed DDs, a specialized branch-and-bound method is designed to successively refine relaxed DDs until proving optimality.

This innovative approach to restructure the graph of the solution set makes DDs a competitive alternative to traditional divide-and-conquer solution techniques such as branch-and-bound and dynamic programming.
Various studies in the past decade have been devoted to showing remarkable improvements of solution time and quality achieved by DDs when compared to the outcome of modern solvers. 
DDs have found their way to a broad array of application areas, from healthcare to supply chain management to finance \citep{bergman2018discrete}.
Other directions of application include cutting plane theory \citep{davarnia:va:2020}, multi-objective optimization \citep{bergman2016multiobjective}, post-optimality analysis \citep{serra:ho:2020}, and integrated branch-and-bound \citep{gonzalez:ci:lo:ro:2020}.

While DDs provide a powerful optimization tool, their applicability domain is limited to discrete problems. 
This limitation is due to a structural requirement that DDs contain a finite number of arcs whose labels represent variable values on their relative domain.
Consequently, a successful extension of DDs to model MIPs has been lacking in the literature.  
\cite{davarnia2020strong} took the initiative in this direction by introducing a new concept called \textit{arc-reduction} that enables building relaxed DDs for continuous nonlinear programs.
Such a DD is used in a cut-generating oracle to obtain linear valid inequalities for the continuous feasible region of the underlying problem.
When embedded inside an outer-approximation scheme, the resulting cutting planes can improve the bounds obtained from classical methods.
This cut-generating framework can be viewed as an interface between DDs and traditional cutting plane methods.
Such an intermediary role, however, does not utilize the full potential of DDs achieved through specialized branch-and-bound and refinement methods.
In the present paper, we develop a framework that \textit{directly} models MIPs through DDs.

\subsection{Background on Unit Commitment Problem} \label{subsec:background_UCP}
The UCP is strongly NP-hard as shown in \cite{bendotti2019complexity}.
Over the past three decades, numerous exact and heuristic approaches have been proposed to solve different variants of the UCP and its extensions under both deterministic and stochastic settings, resulting in a rich literature. 
Examples to solve stochastic/robust UCP include Benders decomposition~\citep{li2007risk}, two-stage stochastic programming~\citep{blanco2017efficient}, heuristic scenario reduction methods~\citep{feng2016solution}, two-stage adaptive robust optimization~\citep{bertsimas2012adaptive}, and algorithms based on constraint generation~\citep{lorca2016multistage}. We encourage the interested reader to consult a survey by~\cite{van2018large} on solution methods for the UCP under uncertainty. Meanwhile, examples to solve deterministic variants include priority listing methods~\citep{lee1992multi}, Tabu search~\citep{rajan2003neural}, Genetic algorithms~\citep{maifeld1996genetic}, dynamic programming~\citep{pang1981evaluation,frangioni2006solving}, simulated annealing~\citep{mantawy1998simulated}, Lagrangian relaxation~\citep{baldick1995generalized}, fuzzy systems~\citep{saneifard1997fuzzy}, and MIP formulations. We refer the interested reader to~\cite{padhy2004unit} for a detailed account on deterministic UCP. 
Other solution approaches include polyhedral analysis of UCP variants; see \cite{rajan2005minimum},~\cite{queyranne2017tight}, and~\cite{bendotti2018min} for examples. 

\cite{garver1962power} is among the pioneers to model the initial variants of the UCP by a MIP that uses three sets of binary variables: (i) on/off variables, (ii) start-up variables, and (iii) shut-down variables. 
Since then, many efforts have been made to improve MIP formulations either by making them easier to solve through reducing the number of variables, constraints and nonzeros, or by proposing stronger formulations with tighter LP relaxation. 
\cite{ostrowski2011tight} propose a new MIP formulation with three sets of binary variables and show that the computational results are superior to those obtained by using fewer sets of 0-1 variables. 
Meanwhile, as explained by the authors, one can generate strong valid inequalities when using all three sets of 0-1 variables, which in turn leads to tighter LP relaxations. 
Subsequently, \cite{morales2013tight} propose a stronger MIP formulation than those of \cite{carrion2006computationally} and \cite{ostrowski2011tight} with a fewer number of constraints and nonzeros. 
Interested reader is referred to~\cite{wu2011tighter} and~\cite{frangioni2008tighter} for other examples of MIP models. 
In addition, a comprehensive analysis and comparison of different MIP formulations are provided by~\cite{knueven2020mixed}.

Due to the substantial cost savings resulted from high quality solutions of the unit commitment,
a host of UCP formulations and solution techniques have been proposed to solve real-world instances; see \cite{bertsimas2012adaptive}, \cite{atakan2017state}, \cite{franz2020long}, and \cite{knueven2020novel}, 
for successful implementations.

\subsection{Contributions} \label{sec:Our_Contributions}

As discussed earlier, the existing DD-based solution methods are limited to discrete problems only.
In this paper, we introduce a framework that generalizes the concept of DDs to solve MIPs directly, which allows for taking advantage of the full arsenal of branch-and-bound and refinement techniques available for DDs. 
In particular, we propose a geometric concept based on decomposing a MIP solution set into rectangular formations of certain properties, which provides both necessary and sufficient conditions for a general optimization problem to be \textit{representable} by DDs. 
The significance of this contribution is two-fold: (i) it notably expands the applicability domain of DDs by lifting the integrality barrier; and (ii) it addresses a fundamental open question in the community on which problem structures are amenable to a DD representation. 
As a consequence, we show that a bounded mixed integer linear program admits a DD representation when reformulated through a specialized Benders decomposition.
This result opens a new vein of applications for DDs to solve a broader class of optimization problems.  
As an evidence for such applications, we develop a novel DD-based solution method to solve the UCP in the energy sector.
Being the first DD-based solution technique for the UCP, our results exhibit a great potential in solving this challenging problem class more efficiently.

The remainder of the paper is organized as follows.
We develop a rectangular decomposition technique for MIPs and establish its connection to constructing DDs with both integer and continuous variables in Section~\ref{sec:DD}.
We transition from the concept of rectangular decomposition to the UCP application through a specialized Benders decomposition technique presented in Section~\ref{sec:Benders}.
Section~\ref{sec:UCP} is devoted to solving the UCP using the developed DD-based framework.
Concluding remarks are given in Section~\ref{sec:conclusion}.

Due to space restrictions, proofs are omitted from the main part of the paper and can be found in Appendix~\ref{sec:proofs}.
Further, a technical comparison of our framework on the integration of Benders decomposition and DDs with those in the literature is given in Appendix~\ref{sec:bd_literature}.

\section{Decision Diagrams Representation} \label{sec:DD}

In this section, we introduce the concept of rectangular decomposition which bridges a general MIP and DDs.
But first, we give a brief overview on the basics of DDs for optimization.
A comprehensive review can be found in \cite{bergman:ci:va:ho:2016}.

\subsection{Overview on Decision Diagrams} \label{sec:background}

Define $N := \{1,\dotsc,n\}$.
We refer to a DD by $\mc{D} = (\mc{U},\mc{A},l(\cdot))$, where $\mc{U}$ denotes the set of nodes and $\mc{A}$ represents the set of arcs, and function $l: \mc{A} \to \Re$ indicates the label of arcs in $\mc{A}$.
The multi-graph induced by $\mc{D}$ is composed of $n$ arc layers $\mc{A}_1, \mc{A}_2, \dotsc, \mc{A}_n$, and $n+1$ node layers $\mc{U}_1,\mc{U}_2, \dotsc, \mc{U}_{n+1}$.
Node layer $\mc{U}_1$ contains a single \textit{root} node $r$, and node layer $\mc{U}_{n+1}$ contains a single \textit{terminal} node $t$.
We define the \textit{width} of a DD as the maximum number of nodes at any node layer,  \blue{that is $\max\{|\mc{U}_i|: i\in\{1,\dots,n+1\}\}$.}

DD $\mc{D}$ represents a set of points of the form $\vc{x} = (x_1,\dotsc,x_n)$ with the following characteristics.
The label $l(a)$ of each arc $a \in \mc{A}_j$, for $j \in N$, represents the value of $x_j$.
Each arc-sequence (path) from $r$ to $t$ encodes a specific value assignment to $\vc{x}$ associated with the labels of the arcs on the path.
The collection of points encoded by all paths of $\mc{D}$ is referred to as the solution set of $\mc{D}$, which is denoted by $\Sol(\mc{D})$. 

Consider a bounded integer set $\mc{P} \subseteq \Z^n$.
Set $\mc{P}$ can be expressed by a DD $\mc{D}$ whose collection of all $r$-$t$ paths encodes the solutions of $\mc{P}$, \textit{i.e.,} $\mc{P} = \Sol(\mc{D})$.
Using this definition, we can model discrete problems with DDs.
Consider a bounded integer program $z^* = \max \left\{f(\vc{x}) \suchthat \vc{x} \in \mc{P}\right\}$
where $f:\Re^n \to \Re$ and $\mc{P} \subseteq \Z^n$. 
To model the above integer program with a DD, we first construct a \textit{weighted} DD, denoted by $[\mc{D}|w(\cdot)]$, where (i) $\mc{D}$ represents an exact DD encoding solutions of $\mc{P}$, and (ii) $w(\cdot):\mc{A} \to \Re$ is a weight function associated with arcs of $\mc{D}$ so that for each $r$-$t$ path $\mr{P} = (a_1,\dotsc,a_n)$, its weight $w(\mr{P}) \coloneqq \sum_{i=1}^{n} w(a_i)$ is equal to $f(\vc{x}^{\mr{P}})$, the objective value of the integral solution corresponding to $\mr{P}$.
Construction of the weight function is immediate when $f(\vc{x})$ is separable, i.e., $f(\vc{x}) = \sum_{i=1}^n f_i(x_i)$, as we can simply define $w(a_i) = f_i(l(a_i))$ for $a \in \mc{A}_i$ and $i \in N$.
The above definition implies that $z^*$ is equal to the weight of a longest path from $r$ to $t$ in the corresponding weighted acyclic graph.

\subsection{Rectangular Decomposition} \label{subsec:PR}

DDs are often described as a union of arc-sequences (paths) from $r$ to $t$.
Since arcs on the path encode certain value assignments to variables, DDs are composed of a finite number of solution points.
The key to extend DDs to model continuous variables is to view them as a union of \textit{node-sequences} from $r$ to $t$. 
From this perspective, multiple arcs with different label values can be considered simultaneously between two consecutive nodes. 
We will show later that, in an extreme case, the arcs between two nodes can virtually span an entire continuous interval between two values in the domain of variables.
As a result, the node-sequence will encode a hyper-rectangle, as opposed to a single point encoded by the arc-sequence.
This analogy lays the foundation for a decomposition technique, which we refer to as \textit{rectangular decomposition}, that restructures the feasible region through a collection of hyper-rectangles with specific properties.
This decomposition framework plays a key role in our ability to represent general sets via DDs.

For any $I \subseteq N$, let $\vc{x}_I$ represent coordinates of $\vc{x} \in \Re^n$ whose index belongs to $I$. 
Given a set $P \subseteq \Re^n$, $\conv(P)$ describes the convex hull of $P$, $\proj_{x_I}(P)$ represents the projection of $P$ onto the space of variables $\vc{x}_I$, and $\dim P$ indicates the dimension of $\conv(P)$. 
Given a function $f(\vc{x}): \Re^n \to \Re$ and a variable subset $I \subseteq N$, we say that $f(\vc{x})$ \textit{is convex in} $\vc{x}_I$ if the restriction of $f(\vc{x})$ to the hyperplane defined by $\{x \in \Re^n | x_i = \bar{x}_i, \, \forall i \in N\setminus I\}$ is convex for any assignment values $\bar{x}_i \in \Re$. 
When we refer to extreme points of a bounded set $P$, denoted by $\mc{X}(P)$, we mean extreme points of $\conv(P)$. 

\begin{theorem} \label{thm:equivalence_general}
Consider a \blue{compact} set $\mc{P} \subseteq \Re^n$, and select $I \subseteq N$.
Assume that there exists a finite collection of \blue{compact} sets $P_I^j$ for $j \in J$, where $J$ is an index set, such that
\begin{itemize}
\item [(i)] $\dim \proj_{x_i}(P_I^j) = 0$, for $i \in N \setminus I$ and $j \in J$, i.e., coordinate $x_i$ is fixed.
\item [(ii)] $\bigcup_{j \in J} \mc{X}(P_I^j) \subseteq \mc{P} \subseteq \bigcup_{j \in J} P_I^j$.
\item [(iii)] For each $j \in J$, there exists a finite collection of hyper-rectangles of the form $R_j^k = \prod_{i=1}^n [l_i^k, u_i^k]$ for $k \in K_j$, where $K_j$ is an index set, such that $\conv(P_I^j) = \conv(\bigcup_{k \in K_j}R_j^k)$.
\end{itemize}
Then, $\max\{f(\vc{x}) | \vc{x} \in \mc{P}\} = \max\{f(\vc{x}) | \vc{x} \in \bigcup_{j \in J} \mc{X}(P_I^j)\} = \max\{f(\vc{x}) | \vc{x} \in \bigcup_{j \in J} \bigcup_{k \in K_j} R_j^k\}$ for any function $f(\vc{x})$ that is convex in $\vc{x}_I$. \Halmos
\end{theorem}

We say that set $\mc{P}$ \textit{admits} a rectangular decomposition w.r.t. $I$, if there exists a finite collection of compact sets $P_I^j$ that satisfy conditions (i)--(iii) of Theorem~\ref{thm:equivalence_general}.
We next illustrate this decomposition concept in an example.

\begin{figure}[!t]
	\centering
	\includegraphics[scale=0.50]{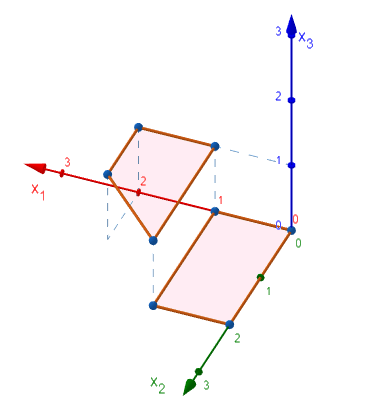}
	\caption{Set $\mc{P}$ of Example~\ref{ex:decompose}}\label{fig:f2}
\end{figure}

\begin{example} \label{ex:decompose}
    Consider a compact mixed integer set 
    \begin{equation*}
        \mc{P} = \left\{(x_1, x_2, x_3) \in [0,2]^2 \times \{0,1\} \suchthat  x_1 + x_2 \leq 3; \,  x_3 - x_1 \leq 0; \, x_1 - x_3 \leq 1  \right\}.
    \end{equation*}
    Set $\mc{P}$ is shown in Figure~\ref{fig:f2}.
    Select $I = \{1, 2\}$.
    To use the result of Theorem~\ref{thm:equivalence_general}, we first introduce sets $P_I^1 = \{\vc{x} \in \Re^3 \suchthat 0 \leq x_1 \leq 1; 0 \leq x_2 \leq 2; x_3 = 0 \}$, and $P_I^2 = \{\vc{x} \in \Re^3 \suchthat 1 \leq x_1 \leq 2; 0 \leq x_2 \leq 2; x_1 + x_2 \leq 3;  x_3 = 1 \}$.
    It is easy to verify that conditions (i) and (ii) of Theorem~\ref{thm:equivalence_general} hold for these sets \blue{as $x_3$ is fixed in sets $P_I^1$ and $P_I^2$, and $\{(0,0,0),(1,0,0),(0,2,0),(1,2,0),(1,1,1),(1,2,1),(2,0,1),(2,1,1)\}\subseteq \mc{P} \subseteq P_I^1 \cup P_I^2$}.
    Next, define \blue{$J=\{1,2\}, K_1=\{1\}, K_2=\{1,2\}$, and} hyper-rectangles $R_1^1 = [0, 1] \times [0,2] \times \{0\}$ for $P_I^1$, and $R_2^1 = \{1\} \times [0, 2] \times \{1\}$ and $R_2^2 = \{2\} \times [0,1] \times \{1\}$ for $P_I^2$.
    These hyper-rectangles satisfy condition (iii) of Theorem~\ref{thm:equivalence_general}.
    As a result, $\max\{f(\vc{x}) \suchthat \vc{x} \in \mc{P}\} = \max\{f(\vc{x}) \suchthat \vc{x} \in \mc{X}(P_I^1) \cup \mc{X}(P_I^2) \} = \max\{f(\vc{x}) \suchthat \vc{x} \in R_1^1 \cup R_2^1 \cup R_2^2\}$ for any function $f(\vc{x})$ that is convex in $(x_1, x_2)$.
\end{example}

In view of Theorem~\ref{thm:equivalence_general}, two levels of decomposition are performed.
First, $\mc{P}$ is decomposed into sets $P^j_I$ through fixing variables not in $I$.
Second, each set $P^j_I$ is decomposed into hyper-rectangles through a convex hull description.
The second level can be viewed as a generalization of the standard extreme point decomposition theorem for maximizing convex functions \citep{rockafellar:1970}, as the hyper-rectangles can be chosen to be the extreme points.
We note, however, that the addition of the first level is necessary to account for the non-convex coordinates of the objective function. 
For instance, consider $\mc{P} = \{(x_1, x_2) \in [-1, 1] \times \{-1, 0, 1\} \suchthat \lvert x_1\rvert \leq \lvert x_2\rvert \}$, and select $I = \{1\}$.
Clearly $\mc{X}(\mc{P}) = \{(-1, -1), (-1, 1), (1, -1), (1, 1)\}$.
It follows that the extreme point decomposition alone does not provide the desired equivalence as $\max \{-x_2^2 \suchthat \vc{x} \in \mc{P}\} \neq \max \{-x_2^2 \suchthat \vc{x} \in \mc{X}(\mc{P})\}$.
Further, we will discuss later that the rectangular volume of the decomposition is key to build efficient DDs as compared to the case for extreme points.
Next, we present a partial converse of Theorem~\ref{thm:equivalence_general} as a sufficient condition for a set to admit rectangular decomposition.

\begin{proposition} \label{prop:converse}
    Consider a \blue{compact} set $\mc{P} \subseteq \Re^n$, and select $I \subseteq N$.
    Assume that there exists a finite set $\mc{Q} \subseteq \Re^n$ such that $\max\{f(\vc{x}) \suchthat \vc{x} \in \mc{P}\} = \max\{f(\vc{x}) \suchthat \vc{x} \in \mc{Q}\}$ for every function $f(\vc{x})$ that is convex in $\vc{x}_I$.
    Then, $\mc{P}$ admits a rectangular decomposition w.r.t. $I$. \Halmos
\end{proposition}

The above rectangular decomposition method has two important consequences in DDs context that will be demonstrated in the next two subsections.

\subsection{Equivalence Class} \label{subsec:equivalence}

In this section, we establish an equivalence relation between DDs that produce the same optimal value for a given objective function.
This result sets the stage for modeling continuous variables via DDs.

\begin{lemma} \label{lem:equivalence}
Consider a \blue{compact} set $\mc{P} \subseteq \Re^n$, and select $I \subseteq N$.
Assume that there exists a finite collection of hyper-rectangles of the form $R_I^j = \prod_{k=1}^n [l_k^j, u_k^j]$ for $j \in J$, where $J$ is an index set, such that 
\begin{itemize}
	\item [(i)] $l_i^j = u_i^j$ for all $i \in N \setminus I$ and $j \in J$, 
	\item [(ii)] $\bigcup_{j \in J} \mc{X}(R_I^j) \subseteq \mc{P} \subseteq \bigcup_{j \in J} R_I^j$.
\end{itemize}
Then, $\max\{f(\vc{x}) | \vc{x} \in \mc{P}\} = \max\{f(\vc{x}) | \vc{x} \in \bigcup_{j \in J} \mc{X}(R_I^j)\} = \max\{f(\vc{x}) | \vc{x} \in \bigcup_{j \in J} R_I^j\}$ for any function $f(\vc{x})$ that is convex in $\vc{x}_I$. \Halmos
\end{lemma}

The difference between the result of Lemma~\ref{lem:equivalence} and that of  Theorem~\ref{thm:equivalence_general} stems from the levels of decomposition involved in the process.
In Lemma~\ref{lem:equivalence}, a set is directly decomposed into rectangular components, whereas in Theorem~\ref{thm:equivalence_general}, a preliminary level of decomposition is performed.
This additional level, through definition of sets $P_I^j$, makes it possible to decompose sets that cannot be directly represented through union of hyper-rectangles.
For instance, it is easy to verify that set $\mc{P}$ of Example~\ref{ex:decompose} cannot be directly decomposed into hyper-rectangles prescribed in Lemma~\ref{lem:equivalence}, whereas it admits rectangular decomposition through the use of Theorem~\ref{thm:equivalence_general}.

In view of Lemma~\ref{lem:equivalence}, we say that the hyper-rectangles $R_I^j$ that satisfy the lemma assumptions form an \textit{equivalence class} w.r.t. $I$, as they share the optimal value of the objective function over any set enveloped between $R_I^j$ and their extreme points.
We next present a similar equivalence class for DD formulations. 
For notational convenience, we use label value $l_a$ for an arc $a \in \mc{A}$ as a shorthand for $l(a)$.
The DDs considered here are not limited by the unique arc-label rule, i.e., they can contain several arcs with equal label values at a layer.
\blue{For uniformity of arguments, we refer to a \textit{virtual} DD as one that contains, between two nodes, a collection of arcs with label values that span a closed continuous interval; see \cite{davarnia2020strong} for an example.}

Consider a DD $\mc{D} = (\mc{U},\mc{A},l(\cdot))$.
For any pair $(u,v)$ of nodes of $\mc{D}$, define $A(u,v)$ to be the set of all arcs of $\mc{D}$ directed from $u$ to $v$.
Further, define $l_{(u,v)}^{\text{max}}$ (resp. $l_{(u,v)}^{\text{min}}$) to be the maximum (resp. minimum) label of the arcs in a non-empty set $A(u,v)$.

\begin{proposition} \label{prop:reduction}
	Consider a DD $\mc{D} = (\mc{U},\mc{A},l(\cdot))$, and select $I \subseteq N$.
	\begin{itemize}
		\item [(i)] Let $\bar{\mc{D}}$ be a DD constructed from $\mc{D}$ with a difference that, for any node pair $(u, v) \in \mc{U}_i \times \mc{U}_{i+1}$ and any $i \in I$, only the arcs with label values $l_{(u,v)}^{\text{min}}$ and $l_{(u,v)}^{\text{max}}$ are maintained and the rest are removed.
		\item [(ii)] Let $\hat{\mc{D}}$ be a virtual DD constructed from $\mc{D}$ with a difference that, for any node pair $(u, v) \in \mc{U}_i \times \mc{U}_{i+1}$ and any $i \in I$, the collection of arcs with label values spanning the interval $[l_{(u,v)}^{\text{min}}, l_{(u,v)}^{\text{max}}]$ is added between $u$ and $v$.
	\end{itemize}
	Then, $\max\{f(\vc{x}) | \vc{x} \in \Sol(\mc{D})\} = \max\{f(\vc{x}) | \vc{x} \in \Sol(\bar{\mc{D}})\} = \max\{f(\vc{x}) | \vc{x} \in \Sol(\hat{\mc{D}})\}$ for any function $f(\vc{x})$ that is convex in $\vc{x}_I$. \Halmos
\end{proposition}

Proposition~\ref{prop:reduction} defines an equivalence class for DDs. 
In words, when there are multiple arcs between two consecutive nodes of a DD, as long as the objective function is convex in the variable corresponding to that arc layer, we can add or remove the middle arcs, and yet preserve the optimal value.
This equivalence property enables us to model continuous intervals by translating virtual DDs into regular DDs. 
This property can also be useful from a computational perspective as the size of the DDs can be reduced through removing unnecessary arcs or through adding arcs that provide more efficient merging possibilities.

\subsection{DD-Representable Sets} \label{subsec:dd_representable}

An important and fundamental question in DDs community concerns the characteristics of optimization problems that can be represented through a DD formulation.

\begin{definition} \label{def:dd_representable}
We say that a \blue{compact} set $\mc{P} \subseteq \Re^n$ is \textit{DD-representable} w.r.t. an index set $I$, if there exists a DD $\mc{D}$ such that $\max\{f(\vc{x}) | \vc{x} \in \mc{P}\} = \max\{f(\vc{x}) | \vc{x} \in \Sol(\mc{D})\}$ for every function $f(\vc{x})$ that is convex in $\vc{x}_I$. 
\end{definition}

In the above definition, the requirement that a set and its associated DD must have the same optimal value for \textit{every} objective function---under appropriate convexity assumptions---is critical for DD-representability.
There are several DD-based procedures, such as Lagrangian relaxation methods \citep{bergman2015lagrangian}, that use DDs with varying objective functions in different stages of the procedure, showing the necessity of properties that guarantee the conservation of optimal values throughout stages.

The question of whether a given set is DD-representable is straightforward when the set is bounded and discrete.
In this case, every point of the set can be encoded by a unique $r$-$t$ path of a DD, leading to a finite width limit.
This argument, however, does not hold when the set contains continuous variables, as all solutions cannot be fully encoded by finitely many paths in a DD.  
In light of the rectangular decomposition technique introduced above, we next give necessary and sufficient conditions for a mixed integer set to be DD-representable.

\begin{corollary} \label{cor:cont_DD}
	Consider a \blue{compact} set $\mc{P} \subseteq \Re^n$, and select $I \subseteq N$.
	The set $\mc{P}$ is DD-representable w.r.t. $I$ if and only if it admits a rectangular decomposition w.r.t. $I$. \Halmos
\end{corollary}

As mentioned above, the conditions of Corollary~\ref{cor:cont_DD} are immediately satisfied for bounded discrete sets.
For mixed integer sets, there is an important class that also satisfies the conditions of Corollary~\ref{cor:cont_DD} as given next.

\begin{corollary} \label{cor:mip}
Let $\mc{Q} \subseteq \Re^{n+1}$ be a compact set and define the mixed integer set $\mc{P} = \{(\vc{x}; y) \in \mc{Q} \, | \, \vc{x} \in \Z^n \}$.
Then, $\mc{P}$ is DD-representable w.r.t. $I = \{n+1\}$. \Halmos
\end{corollary}

The class of mixed integer sets that contain a single continuous variable plays an important role in optimization as it can be viewed as the bridge between a pure discrete case and a general mixed integer case; see e.g. the derivation of mixed integer rounding inequalities that is built upon this analogy \citep{nemhauser1990recursive, conforti2014integer}.
The ability to model this class by DDs, therefore, opens pathways to a wide range of new application domains.
The approach is to represent the continuous components of a general mixed integer set by a new variable and develop its corresponding DD, while accounting for the relation between the continuous components and the substitute variable separately.
The Benders decomposition technique provides an ideal framework for such an approach.
We will establish this framework in Section~\ref{sec:Benders} and apply it to the unit commitment problem in Section~\ref{sec:UCP}. 


In view of Corollary~\ref{cor:mip}, we remark that when a set involves multiple continuous variables, the existence of a rectangular decomposition is not guaranteed.
For example, the unit disk in $\Re^2$ described by $x_1^2 + x_2^2 \leq 1$ does not admit a rectangular decomposition; hence it is not DD-representable.

We conclude this section by presenting a key role of the rectangular decomposition technique in determining the DD size, as the main factor in designing efficient DD-based solution methods.
This implication follows from the fact that the rectangular decomposition of a set is not unique, and it can be achieved in different ways.
For each decomposition, the resulting hyper-rectangles form an equivalence class for DDs representing them.
This variety can be exploited to build DDs with smaller sizes that are computationally more efficient, as illustrated in the next example.

\begin{figure}[!t]
	\centering
	\includegraphics[scale=1.0]{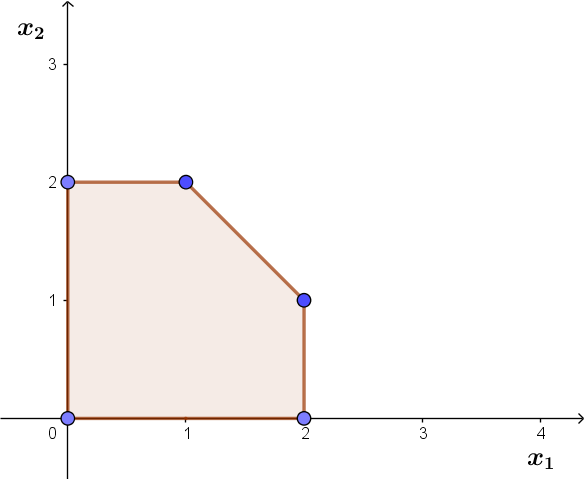}
	\caption{Set $\proj_{(x_1,x_2)}(\mc{P})$ of Example~\ref{ex:equivalence}}\label{fig:f1}
\end{figure}

\begin{figure}[!hbt]
	\centering
	\begin{subfigure}[b]{0.45\linewidth} 
		\centering
		\includegraphics[scale=0.5]{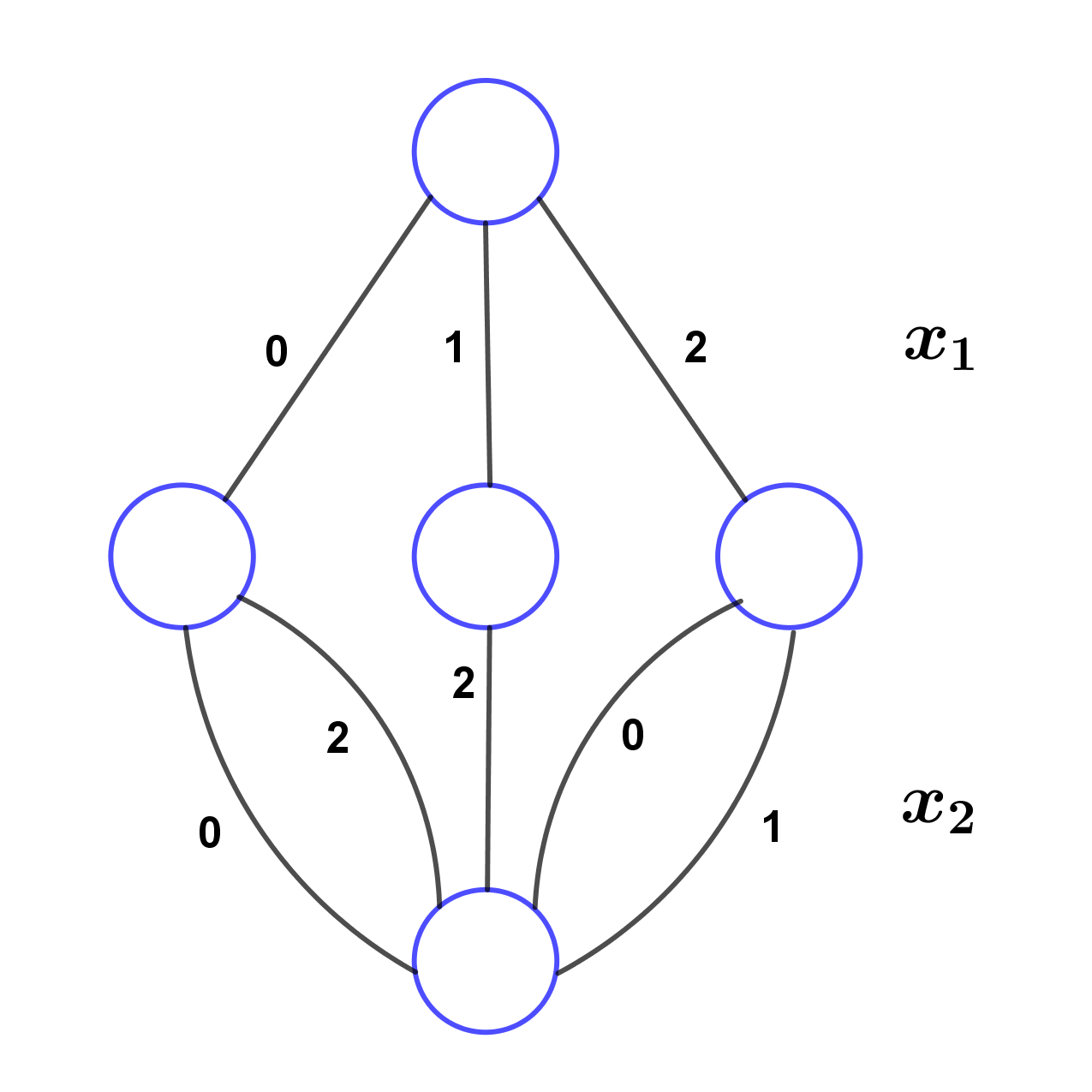}
		\caption{$\mc{D}_1$ with width 3}
		\label{fig:sa1}
	\end{subfigure}
	\hspace{0.05\linewidth}
	\begin{subfigure}[b]{0.45\linewidth} 
		\centering
		\includegraphics[scale=0.5]{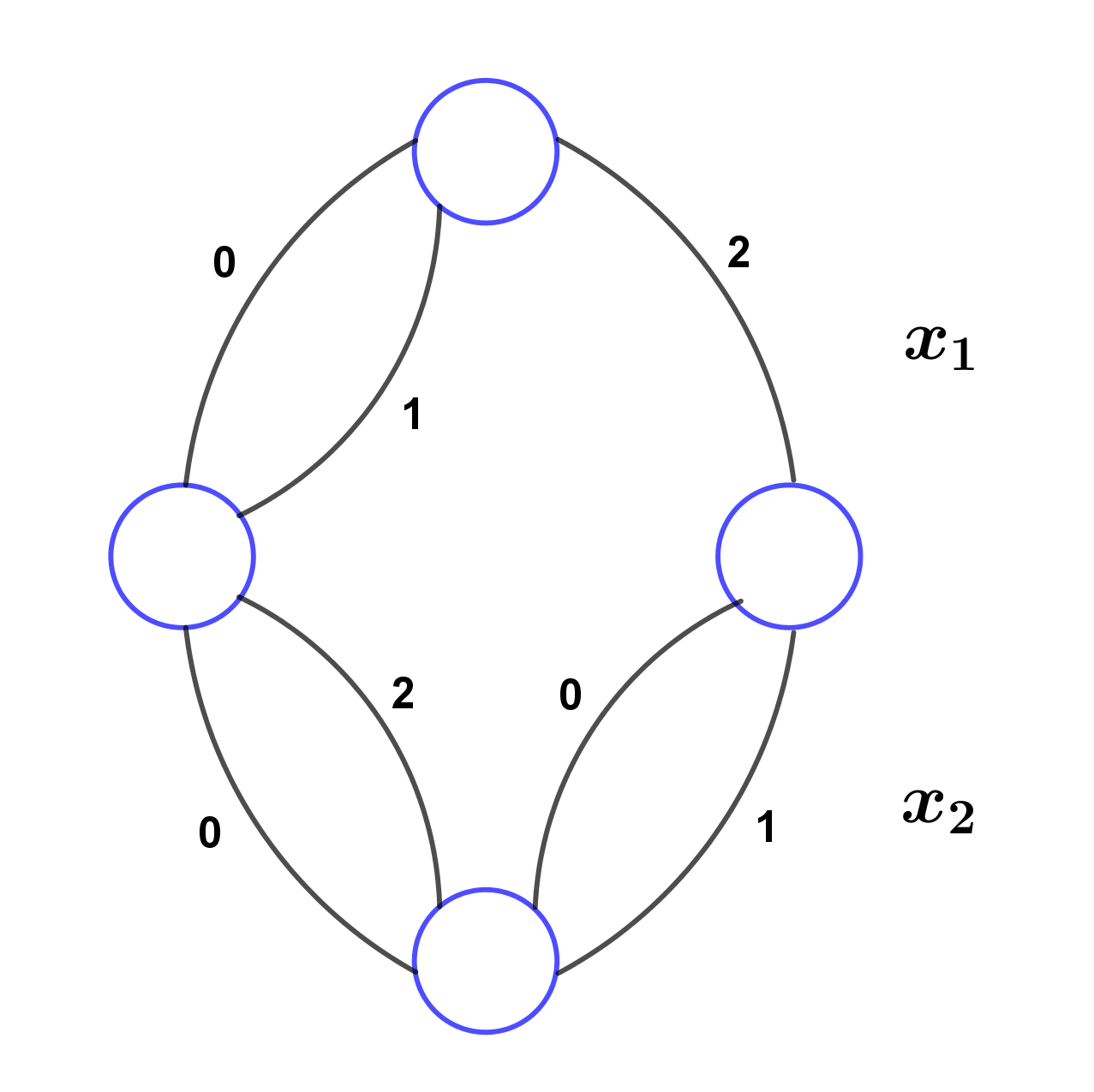}
		\caption{$\mc{D}_2$ with width 2}
		\label{fig:sa2}
	\end{subfigure}	
\caption{The impact of equivalence class on DDs width. Numbers next to arcs represent the labels.}
\label{fig:sa}
\end{figure}

\begin{example} \label{ex:equivalence}
	Consider the set defined by $\proj_{(x_1, x_2)}(\mc{P})$ of Example~\ref{ex:decompose}; see Figure~\ref{fig:f1}.
	The goal is to build a DD $\mc{D}$ such that $\max\{f(\vc{x}) | \vc{x} \in \mc{P}\} = \max\{f(\vc{x}) | \vc{x} \in \Sol(\mc{D})\}$ for every convex function $f(\vc{x})$.
	Define $P_I^1 = \mc{P}$ with $I = \{1,2\}$.
	An immediate rectangular decomposition is obtained through rectangles $R_1^j$ that represent extreme points $\{(0, 0), (0, 2), (1, 2), (2, 0), (2, 1)\}$ of $\mc{P}$.
	The DD $\mc{D}_1$ that models these extreme points forms an equivalence class.
	The minimum width of $\mc{D}_1$ (with the natural ordering of variables in layers) is three, as shown in Figure~\ref{fig:sa1}.
	As an alternative decomposition, we may define hyper-rectangles $R_1^1 = [0,1] \times [0,2]$ and $R_1^2 = \{2\} \times [0,1]$.
	The reduced DD $\mc{D}_2$ representing this equivalence class has width two, Figure~\ref{fig:sa2}.
	This shows the practical advantage of DDs as a member of equivalence classes obtained from different rectangular decompositions.
\end{example}

It follows from Example~\ref{ex:equivalence} that to reduce the size of the DD representation of a given set, one may seek the solution to the combinatorial problem that finds the minimum number of hyper-rectangles that achieve a rectangular decomposition of the set.
This observation implies that a decomposition that contains higher-dimensional hyper-rectangles is generally preferred to one that is simply composed of extreme points (i.e., zero-dimensional rectangles) of the set, as the former can be modeled via fewer node-sequences in a DD.

\section{Benders Decomposition} \label{sec:Benders}

As an application of the rectangular decomposition technique developed in Section~\ref{subsec:PR}, we next present a specialized Benders decomposition (BD) framework that uses DDs to solve MIPs.

Consider a bounded MIP $\mc{H} = \max \{\vc{a}\vc{x} + \vc{b}\vc{y} \, | \, A\vc{x} + B\vc{y} \leq c, \, \vc{x} \in \Z^n \}$ where $\vc{a}$ and $\vc{b}$ are row vectors.
We use BD to define the master problem $\mc{M} = \max \{\vc{a}\vc{x} + z \, | \, (\vc{x}, z) \in \mc{P} \times [l_0, u_0] \}$ where $\mc{P} \subseteq \Z^n$ is the projection of the feasible region of $\mc{H}$ onto $\vc{x}$-space, and the bounds on $z$ are induced from the boundedness of $\mc{H}$.
Subproblems are defined as $\mc{S}(\bar{\vc{x}}) = \max \{ \vc{b}\vc{y} \, | \,  B\vc{y} \leq c - A\bar{\vc{x}} \}$ for a given $\bar{\vc{x}}$.
At each iteration of the algorithm, the output of the subproblems is either a feasibility cut of the form $\vc{\alpha} \vc{x} \leq \alpha_0$ or an optimality cut of the form $z + \vc{\alpha} \vc{x} \leq \alpha_0$, where $\vc{\alpha}$ is a row vector of matching dimension.
These cuts are added to the master problem and the problem is resolved.

It follows from Corollary~\ref{cor:mip} that $\mc{M}$ admits a DD representation, forming an equivalence class.
Such a DD contains $n$ layers corresponding to integer variables $\vc{x}$ and one (last) layer corresponding to the continuous variable $z$ with arc labels representing a lower and an upper bound on this variable. 
\blue{Using this DD, we can find its longest $r$-$t$ path, solve the subproblems for the solution encoding that path, and \textit{refine} the DD with respect to the generated optimality and feasibility cuts.
The refinement operation separates the solutions violating the cuts from the DD solution set; see \cite{bergman:ci:va:ho:2016} for a detailed account on refinement techniques in DD.}
The main difficulty, however, is that the size of a DD that gives an exact representation of $\mc{M}$ grows exponentially with the problem size. 
It is then imperative to employ the notion of \textit{restricted} and \textit{relaxed} DDs in our algorithm to increase efficiency.
The idea is to design DDs with a limited width that represent a restriction (resp. relaxation) of the underlying problem, and thereby providing a primal (resp. dual) bound.
The construction of a restricted DD is straightforward, as it can be achieved by selecting any subset of the $r$-$t$ paths of the exact variant that satisfies the width limit.
The construction of a relaxed DD, however, requires a careful manipulation of the DD structure through a so-called \textit{node-merging} operation in such a way that the resulting DD with a smaller width contains all feasible $r$-$t$ paths of the exact variant, with a possible addition of some infeasible paths.
Such a construction exploits the combinatorial structure of $\mc{M}$, and hence is problem-specific.
We give an instance for designing relaxed DDs for the unit commitment application in Section~\ref{sec:UCP}. 

Algorithm~\ref{alg:DD-BD} presents the full DD-based framework to solve $\mc{H}$ through BD.
In this approach, referred to as DD-BD, we use the following definitions.
Let $C$ be the set of optimality and feasibility cuts, and define $\mc{F}^C$ to be the feasible region described by constraints of $C$.
Consider $\hat{\vc{x}} = (\hat{x}_1, \dotsc, \hat{x}_{|\hat{x}|})$ to be a partial assignment of size $|\hat{x}|$ to variables $x_1, \dotsc, x_{|\hat{x}|}$.
Define $\mc{M}^C(\hat{\vc{x}}) = \max \{\vc{a}\vc{x} + z \, | \, (\vc{x},z) \in \mc{P} \times [l_0, u_0] \cap \mc{F}^C, \, x_i = \hat{x}_i, \forall i = 1,\dotsc,|\hat{x}|\}$ to be the restriction of the master problem $\mc{M}$ after adding the cuts in $C$ and fixing the partial assignment $\hat{\vc{x}}$.
We denote the case with an empty partial assignment and empty constraint set as input by $\mc{M}^{\emptyset}(\emptyset) = \mc{M}$.
We assume that oracles to build relaxed and restricted DDs for the master problem $\mc{M}$ are available.
Using these oracles, one can simply construct relaxed and restricted DDs for $\mc{M}^C(\hat{\vc{x}})$ by fixing the path associated with the partial assignment $\hat{\vc{x}}$ and refining the DD with respect to the cuts in $C$.
Given a DD $\mc{D}$ with $n$ layers, we denote by $\mt{EXACT}(\mc{D})$ the last exact node layer of $\mc{D}$, i.e., the last node layer that contains no merged nodes; see \cite{bergman:ci:va:ho:2016} for properties of the exact cut set operator.

\begin{algorithm}[!htbp]                   		
\caption{DD-BD}
\label{alg:DD-BD}	
\footnotesize
\KwData{$\mc{H}$, $\mc{M}$, $\mc{S}(\bar{\vc{x}})$, and oracles to build relaxed and restricted DDs for $\mc{M}$} 
\KwResult{An optimal solution $(\vc{x}^*, z^*)$ and optimal value $w^*$ of $\mc{H}$} 
initialize set of partial assignments $\hat{\mc{X}}\coloneqq \{\emptyset\}$, set of cuts $C \coloneqq \emptyset$, and lower bound $w^* \coloneqq -\infty$

\While{$\hat{\mc{X}} \neq \emptyset$}{
    pick $\hat{\vc{x}} \in \hat{\mc{X}}$ and update $\hat{\mc{X}} \leftarrow \hat{\mc{X}} \setminus \{\hat{\vc{x}}\}$

    create a restricted DD $\underline{\mc D}$ associated with $\mc{M}^C(\hat{\vc{x}})$
        
    \eIf{$\underline{\mc D} \neq \emptyset$}{
        \Repeat{$\underline{z} = \underline{\rho}$}{
            find a longest $r$-$t$ path in $\underline{\mc D}$ with encoding point $(\underline{\vc{x}}, \underline{z})$ and length $\underline{w}$
                
            solve $\mc{S}(\underline{\vc{x}})$ to obtain feasibility/optimality cuts $\underline{C}$ and optimal value $\underline{\rho}$
                
            update $C \leftarrow C \cup \underline{C}$, and refine $\underline{\mc D}$ w.r.t $\underline{C}$
        }
        \If{$\underline{w} > w^*$}{
            update $w^* \leftarrow \underline{w}$ and $(\vc{x}^*, z^*) \leftarrow (\underline{\vc{x}}, \underline{z})$
        }
    }{
        Go to line 2
    }
 
    create a relaxed DD $\overline{\mc D}$ associated with $\mc{M}^C(\hat{\vc{x}})$
        
    find a longest $r$-$t$ path in $\overline{\mc D}$ with encoding point $(\overline{\vc{x}}, \overline{z})$ and length $\overline{w}$
        
    \eIf{$\overline{w} > w^*$}{
        \Repeat{$\overline{z} = \overline{\rho}$}{
            solve $\mc{S}(\overline{\vc{x}})$ to obtain feasibility/optimality cuts $\overline{C}$ and optimal value $\overline{\rho}$
                
            update $C \leftarrow C \cup \overline{C}$, refine $\overline{\mc D}$ w.r.t $\overline{C}$, and update its longest $r$-$t$ path solution $(\overline{\vc{x}}, \overline{z})$ and length $\overline{w}$
        }
        For each node $u \in \mathtt{EXACT}(\overline{\mc D})$, add the partial assignment encoding a longest $r$-$u$ path of $\overline{\mc D}$ to $\hat{\mc{X}}$
    }{
        Go to line 2
    }
}
Return $(\vc{x}^*, z^*)$ and $w^*$
\end{algorithm}

Next, we show the finiteness and correctness of Algorithm~\ref{alg:DD-BD}, followed by some implementation remarks.

\begin{theorem} \label{thm:DD-BD}
Algorithm~\ref{alg:DD-BD} returns an optimal solution and optimal value of $\mc{H}$ in a finite number of iterations.
\Halmos
\end{theorem}

\begin{remark} \label{rem:DD-BD-1}
For a selected partial assignment $\hat{\vc{x}} \in \hat{\mc{X}}$ in line 3 of Algorithm~\ref{alg:DD-BD}, if the restricted DD $\underline{\mc{D}}$ constructed in line 4 models an exact representation of $\mc{M}^C(\hat{\vc{x}})$, then the construction of a relaxed DD in lines 15--24 can be entirely skipped for that loop.
This follows from the fact that the lower bound at this partial assignment cannot be improved any further by branching down through exact cut sets.
\end{remark}

\begin{remark} \label{rem:DD-BD-2}
The addition of cuts generated from the subproblems to $C$ to be carried over to next iterations is not necessary for the correctness of Algorithm~\ref{alg:DD-BD} as shown in the proof of Theorem~\ref{thm:DD-BD}. 
The advantage of this addition is to avoid generating the same cuts at different iterations for next restricted or relaxed DDs.
In fact, even invoking subproblems for relaxed DDs (in lines 18--21) can be skipped without invalidating correctness of the algorithm.
The addition of this subroutine helps improve the dual bound obtained from the relaxed DDs, and thereby, trigger the pruning condition of line 17 faster.
\end{remark}

The next example illustrates steps of Algorithm~\ref{alg:DD-BD}.

\begin{example}
Consider the following MIP with two integer and two continuous variables.
\begin{align*}
\max_{\vc x\in\{0,1\}^2;~\vc y\in \Re^2_+}\left\{x_1+x_2+2y_1+y_2 \suchthat 
\substack{
x_1 + x_2 \ge 1 \\\\
y_1 + y_2 \ge x_1 + x_2 \\\\
0.3y_1 + 0.7y_2 \le 0.1x_1 + 0.3}\right\}.     
\end{align*}
The optimal solution of this problem is $(x^*_1,x^*_2,y^*_1,y^*_2)=(1,0,1.33,0)$ with the optimal value $3.66$. We intend to solve this problem using the DD-BD approach of Algorithm~\ref{alg:DD-BD}. The master problem $\mc{M}$ is formulated as $\max_{\vc x\in\{0,1\}^2}\left\{x_1 + x_2 + z\suchthat x_1 + x_2 \ge 1\right\}$, where $z$ represents the objective value of the subproblem $\mc{S}(\bar x_1, \bar x_2)$ described by
\begin{align*}
\max_{\vc y\in \Re^2_+}\left\{2y_1 + y_2 \suchthat 
\substack{
y_1 + y_2 \ge \bar x_1 + \bar x_2 \\\\
0.3y_1 + 0.7y_2 \le 0.1\bar x_1 + 0.3
}\right\}.    
\end{align*}
Defining the dual vector $\vc \pi\in \Re^2_+$, we obtain the dual \blue{of the above subproblem as follows}. 
\begin{align*}
\min_{\vc \pi\ge 0}\left\{-(\bar x_1 + \bar x_2)\pi_1 + (0.1 \bar x_1 + 0.3)\pi_2 \suchthat
\substack{
-\pi_1 + 0.3\pi_2 \ge 2 \\\\
-\pi_1 + 0.7\pi_2 \ge 1
}\right\}.
\end{align*}
The feasibility cuts are of the form $\tilde \pi_1(x_1 + x_2)-\tilde \pi_2(0.1x_1 + 0.3) \le 0$ where $\tilde \pi$ is a recession ray of the dual subproblem. Similarly, the optimality cuts are of the form $z + \hat \pi_1(x_1 + x_2)-\hat \pi_2(0.1x_1 + 0.3) \le 0$ where $\hat \pi$ is a feasible solution of the dual problem. We create the exact DD $\underline{\mc D}$ of Figure~\ref{subfig:it1} to represent $\mc{M}$ where $-M$ and $M$ are assumed to be some valid bounds on variable $z$. The longest path of $\underline{\mc D}$, at this iteration, is associated with the solution $(\underline{x_1}, \underline{x_2}, \underline{z})=(1,1,M)$ with length $\underline{w}=2+M$. Solving $S(\underline{\vc x})$ gives a feasibility cut $0.66x_1 + x_2 \le 1$. We update the set of cuts $C$ and refine $\underline{\mc D}$ w.r.t to the feasibility cut. The output is the new DD depicted in Figure~\ref{subfig:it2}. At this iteration, a longest path corresponding to the solution $(\underline{x_1}, \underline{x_2}, \underline{z})=(0,1,M)$ gives the optimal solution $\hat{\vc\pi}=(0,6.66)$ and the optimal value $\underline{\rho}=2$ for the dual subproblem, generating the optimality cut $z \le 0.66x_1 + 2$ through solving $S(\underline{\vc x})$. Refining $\underline{\mc D}$ for the last time w.r.t this optimality cut, the DD of Figure~\ref{subfig:it3} is obtained. It follows that the longest path has length $\underline{w}$ equal to $3.66$ and is achieved by $(\underline{x_1}, \underline{x_2}, \underline{z})=(1,0,2.66)$. At this point, we update $w^*=\underline{w}=3.66$ and $(\vc{x}^*,z^*)=(\underline{\vc x},\underline{z})=(1,0,2.66)$. Since $\underline{\mc D}$ models an exact representation of $\mc M^C(\vc x)$, it follows from Remark~\ref{rem:DD-BD-1} that we can skip creating relaxed DDs, which leads to terminating the algorithm with the optimal solution $(\vc{x}^*,z^*)=(1,0,2.66)$ and the optimal value $w^*=3.66$.  

\begin{figure}
\centering
\begin{subfigure}{0.2\linewidth}
\begin{tikzpicture}[scale=0.7]
\begin{scope}[every node/.style={circle,color=blue,thick,draw,minimum size=0.5cm}]
\node(1) at (0,0) {};
\node(2) at (1.8,-2) {};
\node(3) at (-1.8,-2) {};
\node(4) at (1.8,-4) {};
\node(5) at (-1.8,-4) {};
\node(6) at (0,-6) {};
\end{scope}

\node[text width=2cm, anchor=below, right,scale=1] at (-0.2,-4.2) {\scriptsize $-M$};
\node[text width=2cm, anchor=below, right,scale=1] at (-1,-4.17) {\scriptsize $M$};

\begin{scope}
every edge/.style={draw=red,very thick,}]
\path (1) edge node [right] {\scriptsize 1} (2);
\path (1) edge node [left] {\scriptsize 0} (3);
\path (2) edge[bend left=35] node [right] {\scriptsize 1} (4);
\path (2) edge[bend right=35] node [left] {\scriptsize 0} (4);
\path (3) edge node [left] {\scriptsize 1} (5);
\path (4) edge[bend left=35] node [right] {\scriptsize $M$} (6);
\path (4) edge[bend right=35] node [above] {} (6);
\path (5) edge[bend left=35] node [above] {} (6);
\path (5) edge[bend right=35] node [left] {\scriptsize $-M$} (6);
\end{scope}
\end{tikzpicture}
\caption{Iteration 1}
\label{subfig:it1}
\end{subfigure}\hspace{0.05\linewidth}
\begin{subfigure}{0.2\linewidth}
\centering
\begin{tikzpicture}[scale=0.7]
\begin{scope}[every node/.style={circle,color=blue,thick,draw,minimum size=0.5cm}]
\node(1) at (0,0) {};
\node(2) at (1.8,-2) {};
\node(3) at (-1.8,-2) {};
\node(4) at (1.8,-4) {};
\node(5) at (-1.8,-4) {};
\node(6) at (0,-6) {};
\end{scope}

\node[text width=2cm, anchor=below, right,scale=1] at (-0.2,-4.2) {\scriptsize $-M$};
\node[text width=2cm, anchor=below, right,scale=1] at (-1,-4.17) {\scriptsize $M$};

\begin{scope}
every edge/.style={draw=red,very thick,}]
\path (1) edge node [right] {\scriptsize 1} (2);
\path (1) edge node [left] {\scriptsize 0} (3);
\path (2) edge node [right] {\scriptsize 0} (4);
\path (3) edge node [left] {\scriptsize 1} (5);
\path (4) edge[bend left=35] node [right] {\scriptsize $M$} (6);
\path (4) edge[bend right=35] node [above] {} (6);
\path (5) edge[bend left=35] node [above] {} (6);
\path (5) edge[bend right=35] node [left] {\scriptsize $-M$} (6);
\end{scope}

\end{tikzpicture}
\caption{Iteration 2}
\label{subfig:it2}
\end{subfigure}\hspace{0.05\linewidth}
\begin{subfigure}{0.2\linewidth}
\centering
\begin{tikzpicture}[scale=0.7]
\begin{scope}[every node/.style={circle,color=blue,thick,draw,minimum size=0.5cm}]
\node(1) at (0,0) {};
\node(2) at (1.8,-2) {};
\node(3) at (-1.8,-2) {};
\node(4) at (1.8,-4) {};
\node(5) at (-1.8,-4) {};
\node(6) at (0,-6) {};
\end{scope}

\node[text width=2cm, anchor=below, right,scale=1] at (-0.2,-4.2) {\scriptsize $-M$};
\node[text width=2cm, anchor=below, right,scale=1] at (-1,-4.17) {\scriptsize $2$};

\begin{scope}
every edge/.style={draw=red,very thick,}]
\path (1) edge node [right] {\scriptsize 1} (2);
\path (1) edge node [left] {\scriptsize 0} (3);
\path (2) edge node [right] {\scriptsize 0} (4);
\path (3) edge node [left] {\scriptsize 1} (5);
\path (4) edge[bend left=35] node [right] {\scriptsize $2.66$} (6);
\path (4) edge[bend right=35] node [above] {} (6);
\path (5) edge[bend left=35] node [above] {} (6);
\path (5) edge[bend right=35] node [left] {\scriptsize $-M$} (6);
\end{scope}

\node[text width=2cm,anchor=below, right,scale=1.1] at (3,-1)  {$x_1$};
\node[text width=2cm,anchor=below, right,scale=1.1] at (3,-1)  {$x_1$};

\node[text width=2cm,anchor=below, right,scale=1.1] at (3,-3)  {$x_2$};
\node[text width=2cm,anchor=below, right,scale=1.1] at (3,-3)  {$x_2$};

\node[text width=2cm,anchor=below, right,scale=1.1] at (3.1,-5)  {$z$};
\node[text width=2cm,anchor=below, right,scale=1.1] at (3.1,-5)  {$z$};

\end{tikzpicture}
\caption{Iteration 3}
\label{subfig:it3}
\end{subfigure}
\caption{Different iterations of solving the master problem of Example~\ref{ex:MIP}.}
\label{fig:DDs}
\end{figure}
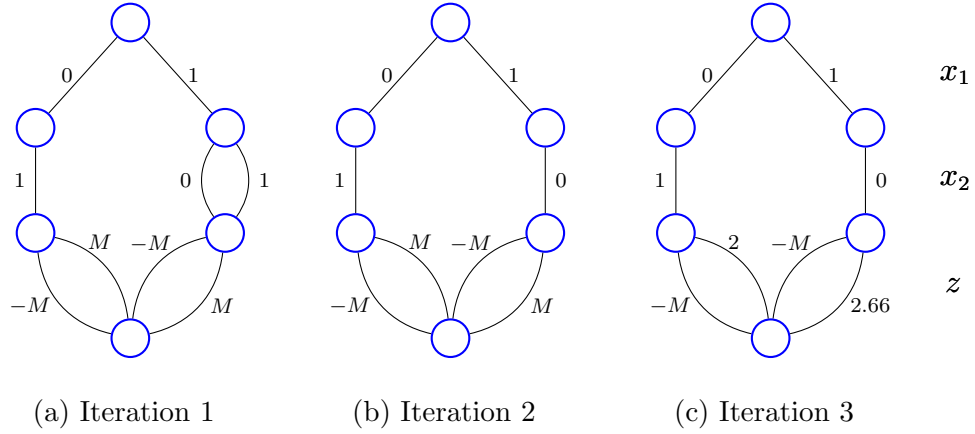

\label{ex:MIP}
\end{example}

While DDs have been used in conjunction with BD for mixed integer models in the literature as in \cite{van2017decision}, their layer construction has been limited to integer variables only, which has led to major computational shortcomings such as not admitting reduced forms---a critical feature for building efficient DDs in practice.
Such limitations are not present in our approach, as it directly incorporates continuous variables within DD layers.
We refer the reader to Appendix~\ref{sec:bd_literature} for a detailed comparison.

\section{Application to the Unit Commitment Problem} \label{sec:UCP}
In this section, we provide an evidence of practicality for the DD-BD framework by applying it to the UCP---a well-known problem in electrid grid applications. 
We first present a common MIP formulation for the UCP.
To apply our framework, we consider the standard BD formulation of the problem, and through a projective transformation, make it amenable to DD representation.
We close the section by comparing the computational results of the DD-BD approach with that of the standard BD formulation.

\subsection{Standard Formulations} \label{subsec:UCP standard}

In the UCP variant we consider in this paper, we seek to find a power generation schedule that minimizes the total operational cost subject to typical UCP constraints such as minimum up/down \blue{times}, generation capacity, demand satisfaction, ramping \blue{requirements}, and spinning reserve.
The typical MIP formulation of this problem is given in \eqref{UCP:opj}--\eqref{const:continuous}; see \cite{rajan2005minimum, ostrowski2011tight, tuffaha2013mixed, bendotti2018min}.
In this formulation, binary variables $x^i_j$ represent whether or not generating unit $i \in N = \{1, \dotsc, n\}$ is up at time $j \in \mc{T} = \{1, \dotsc, T\}$. 
Similarly, binary variables $y^i_j$ (resp. $\bar y^i_j$) indicate whether or not unit $i$ starts up (resp. shuts down) at time $j$.
Further, continuous variables $p^i_j$ and $\bar p^i_j$ denote the production and maximum power available of unit $i$ at time $j$. 

The objective function \eqref{UCP:opj} of the UCP minimizes the total fixed operating costs $c^i_f$, production costs $c^i_g$, and \blue{logarithmic} start-up costs $q^i_j$ for all units over the planning time horizon. 
The fixed operating cost $c^i_f$ is paid whenever unit $i$ is up and working.
The production cost $c^i_g$ of unit $i$ models the variable cost per unit of generated power.
The start-up cost $q^i_j$ of generator $i$ at time $j$ is a \blue{logarithmic} function of the inactive duration of the generator to model the fact that: the colder a generator gets, the greater cost is incurred to start it again. 
This function is commonly linearized by discretizing the time horizon and evaluating the start-up costs $K^i_k$ of unit $i$ after it has been inactive for $k$ consecutive time periods~\citep{ostrowski2011tight, tuffaha2013mixed}.  
This linearization step is modeled in constraint \eqref{const:expo_cost}.
Constraint~\eqref{const:logical} represents the logical relation between the commitment variables $x^i_j$ and variables $y^i_j$ and $\bar{y}^i_j$.
constraints~\eqref{const:min_up} and~\eqref{const:min_down} model the requirement that the schedule of each generator $i$ must satisfy a minimum down-time $\ell^i \geq 1$ and a minimum up-time $L^i \geq 1$, i.e., if a generator $i$ shuts down (resp. starts up) at time $j$, then it must be down (resp. up) for at least $\ell^i$ (resp. $L^i$) time periods. 
In addition, each generator $i$ has a start-up $\SU^i$, ramp-up $\RU^i$, shut-down $\SD^i$, and ramp-down $\RD^i$ rates.
These parameters bound the changes in production rate of generators in consecutive time periods. 
For instance, if unit $i$ is down at time $j$, then its production is bounded above by $\SU^i$ at time $j+1$. 
Similarly, if unit $i$ is working at time $j$ with a generation output strictly greater than $\SD^i$, then it cannot be offline at time $j+1$. 
These requirements are represented in constraints~\eqref{const:ramp_up} and~\eqref{const:ramp_down} as given in \cite{tuffaha2013mixed}.
It is commonly assumed that $\SU^i \leq \RU^i$ and $\SD^i \leq \RD^i$ for all $i \in N$.
Further, each generator $i$ has a minimum $m^i$ and a maximum $M^i$ production capacity that are captured in constraints \eqref{const:pow_gen}.
Finally, constraints \eqref{const:demand} and~\eqref{const:spinning} guarantee that the total demand $D_j$ and the spinning reserve requirement $R_j$ are satisfied at each time period $j$.   
\begin{subequations}
\begin{align}
\min \quad &\sum_{i=1}^n \sum_{j=1}^T c^i_f x^i_j + c^i_g p^i_j + q^i_j \label{UCP:opj}  \\
\text{s.t.} \quad &q^i_j \ge K^i_k\left(x^i_j - \sum_{h=1}^k x^i_{j-h}\right) &\forall k\in\{1,\dots,j-1\}, \quad \forall j\in\mathcal{T}, \quad \forall i\in N \label{const:expo_cost} \\
&y^i_j - \bar{y}^i_j = x^i_j - x^i_{j-1} &\forall j\in\{1,\dots,T\}, \quad \forall i\in N \label{const:logical} \\
&\sum_{j'=j-L^i+1}^j y^i_{j'} \le x^i_j &\forall j\in\{L^i,\dots,T\}, \quad \forall i\in N \label{const:min_up} \\
&\sum_{j'=j-\ell^i + 1}^j \bar{y}^i_{j'} \le 1 - x^i_{j} &\forall j\in\{\ell^i,\dots,T\}, \quad \forall i\in N \label{const:min_down}  \\
&p^i_j - p^i_{j-1} \le \RU^i x^i_{j-1} + \SU^i y^i_j &\forall j\in \mathcal{T}, \quad \forall i\in N \label{const:ramp_up} \\
&p^i_{j-1} - p^i_j \le \RD^i x^i_j + \SD^i \bar{y}^i_j &\forall j\in \mathcal{T}, \quad \forall i\in N \label{const:ramp_down} \\
&m^i x^i_j \le p^i_j \le \bar p^i_j \le M^i x^i_j &\forall j\in \mathcal{T}, \quad \forall i\in N \label{const:pow_gen} \\
&\sum_{i=1}^n p^i_j \ge D_j &\forall j\in \mathcal{T} \label{const:demand}  \\
&\sum_{i=1}^n \bar p^i_j \ge D_j + R_j &\forall j\in \mathcal{T} \label{const:spinning} \\
&x^i_j,y^i_j, \bar{y}^i_j \in\{0,1\} &\forall j\in \mathcal{T}, \quad \forall i\in N \\
&q^i_j,p^i_j,\bar p^i_j \ge 0 &\forall j\in \mathcal{T}, \quad \forall i\in N. \label{const:continuous}
\end{align}
\end{subequations}


The above MIP formulation is used as the core model for the two-stage stochastic UCP, where the first stage decides the commitment status of units, and the second stage determines the generation schedule to satisfy uncertain demand represented by a number of possible scenarios.
Due to the \textit{L-shaped} structure of this stochastic UCP, a BD method is commonly used to solve the problem; see~\citep{zheng2013decomposition, zheng2014stochastic} for examples of such approach.
Other reasons advocating use of BD include slow convergence of the full MIP model because of high memory requirements to store nodes of the branch-and-bound tree and high CPU time needed to solve the linear programming relaxation at each node; see \cite{guan2003optimization, li2005price, fu2013modeling, huang2017electrical} for a detailed exposure.


The standard BD applied to \eqref{UCP:opj}--\eqref{const:continuous} is composed of a master problem that contains binary variables together with the linearized cost variable, and the subproblems that are defined over continuous variables for fixed value assignments to binary variables.
The master problem for this BD formulation is written as
\begin{align}
\min\left\{\sum_{i=1}^n \sum_{j=1}^T \left(c^i_f x^i_j + q^i_j\right) + z  \suchthat \eqref{const:expo_cost}-\eqref{const:min_down}, \,  (\vc{x}, \vc{y}, \vc{\bar y}) \in\{0,1\}^{3nT}, \, \vc{q} \in \Re^{nT}_+, \, z \in [-\Gamma, \Gamma] \right\}, \label{eq:classic BD master}
\end{align}
where $\Gamma$ and $-\Gamma$ are some valid bounds on $z$ induced from the MIP formulation, and $z$ represents the objective value of the following subproblem
\begin{align}
\min\left\{\sum_{i=1}^n\sum_{j=1}^T c^i_g p^i_j \suchthat  \eqref{const:ramp_up}-\eqref{const:spinning}, \, (\vc{p}, \bar{\vc{p}}) \in \Re^{2nT}_+\right\}. \label{eq:classic BD subproblem}
\end{align}
For the two-stage stochastic UCP described above, master problem \eqref{eq:classic BD master} represents the first-stage model, and subproblems \eqref{eq:classic BD subproblem} are defined for each demand scenario realized at the second stage.
\blue{One way to obtain explicit numerical values for $\Gamma$ bounds is to minimize/maximize the objective function of \eqref{eq:classic BD subproblem} over the LP relaxation of \eqref{const:expo_cost}--\eqref{const:continuous}.}
Section~\ref{subsec:Computation} presents computational results for this stochastic problem.

When applying the BD algorithm, at each iteration, an optimal solution of the master problem \eqref{eq:classic BD master} is found and then fed into the subproblem \eqref{eq:classic BD subproblem}, which provides feasibility or optimality cuts to be added back to the master problem.
Since the structure of the above BD formulation conforms to that used in the DD-BD approach, we adopt this model as the basis for our analysis.
We show in the next section that through exploiting the DD structure, we can streamline the master problem representation, which results in a superior computational performance. 

\subsection{DD-BD: Master Problem Formulation} \label{subsec:DDBD UCP master}

The first step in applying the DD-BD approach of Section~\ref{sec:Benders} to the UCP is to construct a DD that represents the feasible region of the master problem \eqref{eq:classic BD master}.
Since this set is defined over binary variables $(\vc{x}, \vc{y}, \bar{\vc{y}})$ together with continuous variables $(\vc{q}, z)$, a typical DD would contain all five variable types in its arc layers.
It turns out that it suffices to model variables $\vc{x}$ and $z$ only in DD arc layers, while recording the value of other variables through state representations.
This projection will allow for a direct application of Algorithm~\ref{alg:DD-BD}, and it will lead to a significant size reduction for DD models.

We present the DD structure for the single-unit case, i.e., $n = 1$.
The extension to multi-unit case follows similarly by replicating the DD structure for other units, as they can be considered independently in the master problem.
For this reason and to simplify notation, we drop the superscript representing unit number in variables and parameters in the sequel. 

The construction of the DD representing \eqref{eq:classic BD master} is given in Algorithm~\ref{alg:DD master}.
In view of this algorithm, each node $u \in \mc{U}_j$ for $j \in \{1, \dotsc, T+1\}$ contains a state value of the form $(s_u^+, s_u^-)$ where $s_u^+$ and $s_u^-$ record the number of time periods passed since the last start-up and shut-down of the unit at time $j$, respectively.
An advantage of this state definition is its \textit{memoryless} property that allows for a direct evaluation of the minimum up/down requirements as well as the \blue{logarithmic} start-up cost at each time period without the need for backtracking to identify the unit status at previous periods; see the proof of Theorem~\ref{thm:DD master equivalence} for a detailed analogy.
The initial state value at the root node is set to $(\infty, \infty)$ to indicate that the unit is down at the start of the planning timeline and is ready to start up if decided to.
The output of Algorithm~\ref{alg:DD master} is a DD $\mc{D} = (\mc{U}, \mc{A}, l(.))$ with $T + 1$ arc layers representing variables $x_j$ for $j \in \mc{T}$, as well as the continuous variable $z$ at the last layer.
Each arc $a \in \mc{A}$ has a label $l_a$ that represents the value assignment of its corresponding variable, and a weight $w_a$ that captures the objective function rate of that assignment. \blue{The following example illustrates an exact DD for the master problem of an instance of the UCP.}

\begin{example}\label{ex:DD-UCP}
\blue{Consider a UCP instance with one generator over two time periods with the fixed cost $c_f$, and the logarithmic start-up cost $K_k$ for $k$ consecutive inactive time periods. Further, suppose that the minimum down-time of the generator is one and the minimum up-time is two. The DD depicted in Figure~\ref{fig:DD-UCP} represents the feasible region of the associated master problem where $-\Gamma$ and $\Gamma$ are valid bounds for the continuous variable $z$. In this DD, the dashed and solid arcs represent arc labels with values of zero and one, respectively. In addition, the state value of each node and weight of each arc are given next to them, where $K_{\infty}$ denotes the cost for the first time that the unit starts up.}

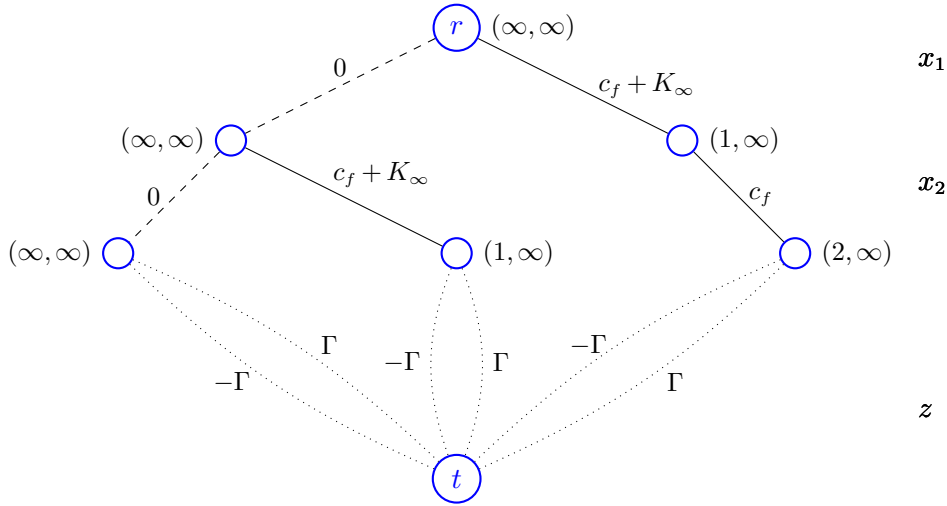
\begin{figure}[!htbp]
\centering
\begin{tikzpicture}[scale=1.5]
\begin{scope}[sn/.style={circle,color=blue,thick,draw,minimum size=0.4cm}]
\node[sn, label=right:{\small $(\infty, \infty)$}](1) at (0,0) {$r$};
\node[sn, label=right:{\small $(1, \infty)$}](2) at (2,-1) {};
\node[sn, label=left:{\small $(\infty, \infty)$}](3) at (-2,-1) {};
\node[sn, label=right:{\small $(2, \infty)$}](4) at (3,-2) {};
\node[sn, label=right:{\small $(1, \infty)$}](5) at (0,-2) {};
\node[sn, label=left:{\small $(\infty, \infty)$}](6) at (-3,-2) {};
\node[sn](7) at (0,-4) {$t$};
\end{scope}
\node[text width=2cm, anchor=below, right,scale=1] at (4,-0.3) {$x_1$};
\node[text width=2cm, anchor=below, right,scale=1] at (4,-0.3) {$x_1$};
\node[text width=2cm, anchor=below, right,scale=1] at (4,-0.3) {$x_1$};
\node[text width=2cm, anchor=below, right,scale=1] at (4,-1.4) {$x_2$};
\node[text width=2cm, anchor=below, right,scale=1] at (4,-1.4) {$x_2$};
\node[text width=2cm, anchor=below, right,scale=1] at (4,-1.4) {$x_2$};
\node[text width=2cm, anchor=below, right,scale=1] at (4,-3.4) {$z$};
\node[text width=2cm, anchor=below, right,scale=1] at (4,-3.4) {$z$};
\node[text width=2cm, anchor=below, right,scale=1] at (4,-3.4) {$z$};
\path(1) edge node [right] {\small $\hspace{2mm}c_f+K_\infty$} (2);
\path[dashed] (1) edge node [above] {\small $0$} (3);
\path(2) edge node [right] {\small $c_f$} (4);
\path(3) edge node [above] {\small $\hspace{10mm}c_f+K_\infty$} (5);
\path(3)[dashed] edge node [left] {\small $\hspace{-5mm}0$} (6);
\path[dotted](4) edge[bend left = 10] node [right] {\small $\hspace{2mm}\Gamma$} (7);
\path[dotted](4) edge[bend right = 10] node [left] {\small $\hspace{-2mm}-\Gamma$} (7);
\path[dotted](5) edge[bend right = 20] node [left] {\small $-\Gamma$} (7);
\path[dotted](5) edge[bend left = 20] node [right] {\small $\Gamma$} (7);
\path[dotted](6) edge[bend left = 10] node [right] {\small $\hspace{2mm}\Gamma$} (7);
\path[dotted](6) edge[bend right = 10] node [left] {\small $\hspace{-10mm}-\Gamma$} (7);
\end{tikzpicture}
\caption{\blue{Illustration of an exact DD in Example~\ref{ex:DD-UCP}}}
\label{fig:DD-UCP}
\end{figure}    
\end{example}

\begin{algorithm}[!ht]                    				  	
	\caption{The construction of DD for the master problem of the UCP for $n=1$}          
	\label{alg:DD master}			                        
		\KwData{parameters $\ell$, $L$, $\Gamma$, $K$}
		\KwResult{a weighted DD $[\mc{D}|w(.)]$}
		create the root node with state value $(\infty, \infty)$.
		
		\ForAll{$j \in \mc{T}$ and $u \in \mc{U}_j$}{
    		\eIf{$s_u^+ \geq s_u^-$}{
    		    create a node $v \in \mc{U}_{j+1}$ with state value $(s_u^+ +1 , s_u^- +1)$, and an arc $a \in \mc{A}_j$ connecting $u$ to $v$ with $l_a = 0$ and $w_a = 0$
    		    
    		    \If{$s_u^- \geq \ell$}{
    		    create a node $v \in \mc{U}_{j+1}$ with state value $(1 , s_u^- +1)$, and an arc $a \in \mc{A}_j$ connecting $u$ to $v$ with $l_a = 1$ and $w_a = c_f + K_{s_u^-}$}
    		}{
    		    create a node $v \in \mc{U}_{j+1}$ with state value $(s_u^+ +1, s_u^- +1)$, and an arc $a \in \mc{A}_j$ connecting $u$ to $v$ with $l_a = 1$ and $w_a = c_f$
    		    
    		    \If{$s_u^+ \geq L$}{
        		create a node $v \in \mc{U}_{j+1}$ with state value $(s_u^+ +1, 1)$, and an arc $a \in \mc{A}_j$ connecting $u$ to $v$ with $l_a = 0$ and $w_a = 0$}
    		}
		}
		\ForAll{$u \in \mc{U}_{T+1}$}{
		    create two arcs $a_1, a_2 \in \mc{A}_{T+1}$ connecting $u$ to the terminal node with $l_{a_1} = w_{a_1} = \Gamma$, $l_{a_2} = w_{a_2}  = -\Gamma$.
		}
\end{algorithm}

The next result shows that the solution set of the equivalence class formed by the DD obtained from Algorithm~\ref{alg:DD master} matches the projection of the feasible region of the master problem~\eqref{eq:classic BD master}.
In what follows, we refer to a solution $(\dot{\vc{x}}, \dot{\vc{y}}, \dot{\bar{\vc{y}}}, \dot{\vc{q}}, \dot{z})$ of \eqref{eq:classic BD master} as \textit{non-redundant} if there is no distinct point $(\dot{\vc{x}}, \dot{\vc{y}}, \dot{\bar{\vc{y}}}, \tilde{\vc{q}}, \dot{z})$ of \eqref{eq:classic BD master} with a strictly smaller objective value.

\begin{theorem} \label{thm:DD master equivalence}
Let $\mc{D}$ be a DD constructed by Algorithm~\ref{alg:DD master} in relation to a UCP with the master problem~\eqref{eq:classic BD master}.
Then, a point $(\dot{\vc{x}}, \dot{z})$ encodes an $r$-$t$  path of length $\dot{c}$ in the equivalence class formed by $\mc{D}$, if and only if this point can be extended to a non-redundant solution $(\dot{\vc{x}}, \dot{\vc{y}}, \dot{\bar{\vc{y}}}, \dot{\vc{q}}, \dot{z})$ of \eqref{eq:classic BD master} with objective value $\dot{c}$. 
\Halmos
\end{theorem}

It follows from Theorem~\ref{thm:DD master equivalence} that refining the DD constructed by Algorithm~\ref{alg:DD master} with respect to optimality and feasibility cuts produced from subproblems~\eqref{eq:classic BD subproblem} results in a solution set equivalent to that of the projection of the non-redundant feasible region of \eqref{eq:classic BD master} after addition of those cuts.
As a consequence, the proposed DD-BD approach converges to an optimal solution of the UCP problem.
We will give details on the derivation of optimality and feasibility cuts in the next section.

Algorithm~\ref{alg:DD master} demonstrates the construction of an \textit{exact} DD representing the master problem \eqref{eq:classic BD master}, which can be used as the DD oracle required for Algorithm~\ref{alg:DD-BD}.
For practical applications, however, the construction of such a DD could be computationally prohibitive.
To mitigate this computational burden, as discussed in Section~\ref{sec:Benders}, it is necessary to use restricted and relaxed DDs with smaller width limits to improve efficiency of the algorithm.
We next give details on designing a relaxed DD for \eqref{eq:classic BD master}.

At the heart of a relaxed DD lies a \textit{node-merging} operation, through which multiple nodes at each layer of the DD are merged into one node, and thereby reducing the DD size.
This node-merging operation must satisfy the condition that any $r$-$t$ path of the exact DD remains an $r$-$t$ path of the resulting DD with no greater length (for a minimization original problem.)
As a byproduct, such a merging operation could create some infeasible $r$-$t$ paths, hence providing a relaxation.
The following definitions set the stage for building a relaxed DD for the UCP.

\begin{definition} \label{def:modified alg}
Consider the following modification of Algorithm~\ref{alg:DD master}: (i) each node $u \in \mc{U}_j$ for $j \in \{1, \dotsc, T+1\}$ is assigned with state values $(s^+_u, s^-_u, s^=_u)$ with an additional component $s^=_u$ whose value is initialized and updated similarly to that of $s^-_u$ throughout the iterations of the algorithm; (ii) parameter $K_{s^-_u}$ is replaced with $K_{s^=_u}$ when calculating arc weights in line 6 of the algorithm.
\end{definition}

\begin{definition} \label{def:merging}
Consider the following node-merging operation defined over the DDs obtained from the algorithm of Definition~\ref{def:modified alg}. At layer $j \in \{1, \dotsc, T+1\}$, a collection of nodes $\{u_1, \dotsc, u_k\}$ with the property that either $s^+_{u_i} \geq s^-_{u_i}$ for all $i \in \{1,\dotsc,k\}$, or $s^+_{u_i} < s^-_{u_i}$ for all $i \in \{1,\dotsc,k\}$ are merged into node $v$ with state values
$s^+_v = \max_{i \in \{1,\dotsc,k\}}s^+_{u_i}$, $s^-_v = \max_{i \in \{1,\dotsc,k\}}s^-_{u_i}$, and $s^=_v = \min_{i \in \{1,\dotsc,k\}}s^=_{u_i}$.
\end{definition} 

We next show that the modified algorithm of Definition~\ref{def:modified alg} combined with the node-merging operation of Definition~\ref{def:merging} yield a relaxed DD for the master problem \eqref{eq:classic BD master}.

\begin{proposition} \label{prop:relaxed DD UCP}
Let $\bar{\mc{D}}$ be the DD obtained from Algorithm~\ref{alg:DD master} in relation to the master problem \eqref{eq:classic BD master}. 
Let $\tilde{\mc{D}}$ be a DD constructed from the modified algorithm of Definition~\ref{def:modified alg} in combination with the node-merging operation of Definition~\ref{def:merging}.
Then, $\tilde{\mc{D}}$ provides a relaxation of $\bar{\mc{D}}$.
\Halmos
\end{proposition}

\subsection{DD-BD: Subproblem Formulation} \label{subsec:DDBD UCP subproblem}

Through application of Algorithm~\ref{alg:DD-BD} to the UCP, subproblem \eqref{eq:classic BD subproblem} is solved for fixed values of master variables to obtain feasibility and optimality cuts, with respect to which the restricted and relaxed DDs representing the master problem are refined.
These cuts, however, are generated in the space of variable $(\vc{x}, \vc{y}, \bar{\vc{y}}, z)$.
As a result, they cannot be directly used to refine the DDs representing the master problem as they only contain variables $(\vc{x}, z)$; see Section~\ref{subsec:DDBD UCP master}.
To resolve this discrepancy, we next show that constraints \eqref{const:ramp_up} and \eqref{const:ramp_down} in the description of the UCP can be replaced with two constraints that do not contain variables $(\vc{y}, \bar{\vc{y}})$.
Such a substitution in the subproblem \eqref{eq:classic BD subproblem} leads to feasibility and optimality cuts that are in the space of $(\vc{x}, z)$, which can be directly used to refine DDs of the master problem.

\begin{proposition} \label{prop:ramp_equivalent}
Assume that $\SU^i \le \RU^i$ and $\SD^i \le \RD^i$ for all generators $i\in N$. 
Let $\mc{E}$ be the feasible region of the UCP described by \eqref{const:expo_cost}--\eqref{const:continuous}.
Define $\mc{G}$ to be the feasible region of the UCP where constraints \eqref{const:ramp_up} and \eqref{const:ramp_down} are respectively replaced with
\begin{subequations}
\begin{align}
&p^i_j - p^i_{j-1} \le \left(\RU^i-\SU^i\right)x^i_{j-1} + \SU^i x^i_j &\forall j\in \mathcal{T}, \quad \forall i\in N \label{eq:ramp_up_2} \\
&p^i_{j-1} - p^i_j \le \left(\RD^i - \SD^i\right)x^i_j + \SD^i x^i_{j-1} &\forall j\in \mathcal{T}, \quad \forall i\in N. \label{eq:ramp_down_2}
\end{align}
\end{subequations}
Then, $\mc{E} = \mc{G}$. \Halmos
\end{proposition}

Substituting \eqref{const:ramp_up} and \eqref{const:ramp_down} with \eqref{eq:ramp_up_2} and \eqref{eq:ramp_down_2} in the description of the subproblem \eqref{eq:classic BD subproblem}, we obtain the following dual problem 
\begin{subequations}
\begin{align}
\max \quad &\sum_{j=1}^T D_j \psi_j + \sum_{j=1}^T\left(D_j + R_j\right)\beta_j  + \sum_{i=1}^n \sum_{j=1}^T  \bar{x}^i_j\left(m^i\phi_{i,j} - M^i \pi_{i,j}\right) \\  
&+ \sum_{i=1}^n \sum_{j=1}^T \Big(\left(\SU^i - \RU^i\right) \bar{x}^i_{j-1} - \SU^i \bar{x}^i_j\Big)\gamma_{i,j} \nonumber \\
&+ \sum_{i=1}^n \sum_{j=1}^T \Big(\left(\SD^i - \RD^i\right) \bar{x}^i_j-\SD^i \bar{x}^i_{j-1}\Big)\delta_{i,j} \nonumber \\
\text{s.t.} \quad &\psi_j - \gamma_{i,j}+\gamma_{i,j+1} + \delta_{i,j}-\delta_{i,j+1} + \phi_{i,j} - \eta_{i,j} \le c^i_g &\forall j\in \mathcal{T}, \quad \forall i\in N \\
&\beta_j + \eta_{i,j} - \pi_{i,j} \le 0 &\forall j\in \mathcal{T}, \quad \forall i\in N \\
&\psi_j, \beta_j, \gamma_{i,j}, \delta_{i,j}, \phi_{i,j}, \eta_{i,j}, \pi_{i,j} \ge 0 &\forall j\in \mathcal{T}, \quad \forall i\in N. 
\end{align}
\end{subequations}

When applying Algorithm~\ref{alg:DD-BD}, we solve the above dual problem for each given assignment $\bar{\vc{x}}$.
If the dual subproblem is unbounded (i.e., the subproblem is infeasible,) we obtain a ray $(\vc{\tilde \psi}, \vc{\tilde \beta}, \vc{\tilde \phi}, \vc{\tilde \pi}, \vc{\tilde \gamma}, \vc{\tilde \delta})$ and add the following feasibility cut to the master problem.
\begin{align}
0 &\ge \sum_{j=1}^T D_j \tilde \psi_j + \sum_{j=1}^T\left(D_j + R_j\right)\tilde \beta_j  + \sum_{i=1}^n \sum_{j=1}^T \left(m^i\tilde \phi_{i,j} - M^i \tilde \pi_{i,j}\right)x^i_j \label{eq:feasibility cut} \\  
&+ \sum_{i=1}^n \sum_{j=1}^T \bigg(\tilde \gamma_{i,j}\Big(\left(\SU^i - \RU^i\right) x^i_{j-1} - \SU^i x^i_j\Big) + \tilde \delta_{i,j} \Big(\left(\SD^i - \RD^i\right) x^i_j-\SD^i x^i_{j-1}\Big)\bigg). \nonumber 
\end{align}

Otherwise, we obtain an optimal solution $(\vc{\psi^*}, \vc{\beta^*}, \vc{\phi^*}, \vc{\pi^*}, \vc{\gamma^*}, \vc{\delta^*})$ of the dual subproblem and add the following optimality cut to the master problem. 
\begin{align}
z &\ge \sum_{j=1}^T D_j \psi^*_j + \sum_{j=1}^T\left(D_j + R_j\right)\beta^*_j + 
\sum_{i=1}^n \sum_{j=1}^T \left(m^i\phi^*_{i,j} - M^i \pi^*_{i,j}\right) x^i_j \label{eq:optimality cut} \\ 
&+ \sum_{i=1}^n \sum_{j=1}^T \bigg(\gamma^*_{i,j} \Big(\left(\SU^i - \RU^i\right) x^i_{t-j} - \SU^i x^i_j\Big) + \delta^*_{i,j} \Big(\left(\SD^i - \RD^i\right) x^i_j-\SD^i x^i_{j-1}\Big)\bigg). \nonumber 
\end{align}

\subsection{Computational Experiments} \label{subsec:Computation}
In this section, we assess the performance of the proposed DD-BD algorithm in solving the stochastic UCP instances for a short-term (day-ahead) power generation model. 
\blue{For our experiments, we consider benchmark instances from the problem class OR-LIB/UC in UnitCommitment.jl~\citep{UnitCommitmentjl}.
Since these instances lack some of the modeling elements in the UCP variant we consider in this paper, we set their associated parameter values as follows. 
We set the production cost $c^i_g$ for unit $i\in N$ as the average of production costs for that unit. 
We create the fixed operating cost $c^i_f$ for each unit $i\in N$ randomly from the interval $[400,1000]$. 
We set $RU^i=SU^i$ and $RD^i=SD^i$ for each unit $i\in N$.
We capture demand uncertainty by a set of scenarios $\Xi$ that are generated randomly from the interval $[0.75 \bar{M}, \bar{M}]$ where $\bar M$ denotes the total capacity over all generators.
}
All experiments were conducted on a machine running Windows 10, x64 operating system with Intel\textsuperscript{\textregistered} Core i7 processor (2.60 GHz) and 32 GB RAM. The code is written in Microsoft Visual Studio 2015 in C++.  
For the optimization parts, we use the Gurobi solver version 9.0.0.

\subsubsection{Comparison with Decomposition Methods} \label{subsubsec:decomposition}
\blue{We first present computational results that compare the outcome of the DD-BD method, as a decomposition technique, with other standard decomposition methods.
In particular, we consider two BD approaches referred to as G-BD 1 and G-BD 2.}

\blue{G-BD 1 employs the classical BD in a \textit{textbook} fashion, where the master problem \eqref{eq:classic BD master} is solved by Gurobi at its default settings, and the optimality/feasibility cuts are produced by solving the subproblems \eqref{eq:classic BD subproblem} for each demand scenario.}

\blue{Since the traditional implementation of BD, as used in G-BD 1, is known to suffer from a slow convergence, we employ modern modeling and algorithmic practices in the literature that boost the performance of the G-BD.
We refer to this improved approach as G-BD 2.
The details of these improvements are as follows. 
(i) To accelerate the convergence rate associated with the feasibility cuts, we use the $L^1$ normalization technique proposed by by~\cite{fischetti2010note} and \cite{bonami2020implementing}. 
Consider the MIP $\mc{H}$ defined in Section~\ref{sec:Benders}. 
Assume that the subproblem $\mc{S}(\bar{\vc{x}})$ is infeasible for a fixed solution $\bar{\vc{x}}$ obtained from the master problem.
It follows that $w^* = \min_{w,\vc{y}}\{w \suchthat B\vc{y} + w\vc{1} \leq \vc{c} - A\bar{\vc{x}}\} = \vc{\pi}^*(\vc{c}-A\bar{\vc{x}}) > 0$, where $\vc{\pi}^*$ is the optimal dual vector of the \textit{normalized} subproblem, $\vc{1}$ is a unit column vector, and the equality holds because of the strong duality property. 
Therefore, the inequality $\vc{\pi}^*(\vc{c}-A\bar{\vc{x}}) \le 0$ can be added to the master problem as a feasibility cut that excludes $\bar{\vc{x}}$.
(ii) During presolve, we generate heuristic cuts to tighten the initial LP relaxation of the master problem as follows. 
At each time period, we identify the scenario with the highest demand and add valid inequalities in the space of commitment variables to ensure the feasibility of constraints~\eqref{const:pow_gen} and~\eqref{const:spinning}. 
(iii) We derive the so-called \textit{combinatorial cuts} that exclude the current infeasible solution of the master problem by forcing a change of value in at least one commitment variable; see~\cite{codato2006combinatorial} and~\cite{rodriguez2021accelerating} for details about these cuts. 
(iv) Finally, we strengthen the logarithmic start-up cost constraints~\eqref{const:expo_cost} by reducing the coefficient of variables $x^i_{j-h}$ as proposed in~\cite{silbernagl2015improving} and~\cite{knueven2020novel}.
}

\blue{For the DD-BD method, we set the width limit for the restricted/relaxed DDs to two and use the node-merging operation given in Definition~\ref{def:merging} for the construction of relaxed DDs. 
When necessary during the refinement procedures, we allow the increase in the DD width to accommodate for the separation of current non-optimal solutions. 
}

\blue{
For these experiments, we fix the number of demand scenarios at 100.
We consider four size categories with $n \in \{20, 50, 100, 150\}$ for the number of generators, and solve five 5 instances for each problem size as given in the OR-LIB/UC.
Table~\ref{tab:1} contains the result of solving these instances with G-BD 1, G-BD 2 and DD-BD. 
For these runs, we fix the time limit at 7200 seconds.
The first column of Table~\ref{tab:1} indicates the number of generators.
The second column shows the instance number.
The columns under ``G-BD 1" report the solution time (in seconds), the total number of generated feasibility and optimality cuts during the solution process.
The subsequent columns contain the same quantities for G-BD 2 and DD-BD.
As observed in Table~\ref{tab:1}, while the boosting methods implemented in G-BD 2 have substantially improved the solution time over G-BD 1, the DD-BD approach uniformly outperforms both G-BD alternatives across all size categories, rendering it as the most effective decomposition approach. 
}

\begin{table}[!htbp]
\begin{center}
\caption{\blue{Solution times (in seconds) of G-BD 1, G-BD 2, and DD-BD for instances of the stochastic UCP with $|\Xi|=100$.}}
\label{tab:1}
\scalebox{1}{
\begin{tabular}{lc|rrr|rrr|rrr} 
\multicolumn{2}{c|}{} & \multicolumn{3}{c|}{{G-BD 1}} & \multicolumn{3}{c|}{{G-BD 2}} & \multicolumn{3}{c}{{DD-BD}}\\
$n$ & Instance & time & f cuts & o cuts & time & f cuts & o cuts & time & f cuts & o cuts\\\hline
20 & 1 & 753.21 & 633 & 34 & 527.76 & 393 & 31 & \textbf{381.78} & 301 & 26 \\
20 & 2 & 981.60 & 638 & 26 & 508.32 & 438 & 28 & \textbf{393.94} & 472 & 34 \\
20 & 3 & 854.45 & 920 & 20 & 499.68 & 586 & 19 & \textbf{371.70} & 560 & 24 \\
20 & 4 & 1516.88 & 196 & 31 & 732.96 & 162 & 29 & \textbf{462.88} & 183 & 21 \\
20 & 5 & 1240.28 & 264 & 34 & 805.59 & 178 & 33 & \textbf{518.86} & 121 & 25 \\
\hline
50 & 1 & 4633.54 & 1119 & 73 & 2683.69 & 939 & 65 & \textbf{1238.59} & 663 & 79 \\
50 & 2 & 4094.20 & 973 & 60 & 2154.84 & 850 & 54 & \textbf{1318.84} & 731 & 40 \\
50 & 3 & 3791.37 & 1428 & 33 & 2071.47 & 899 & 32 & \textbf{1149.24} & 836 & 29 \\
50 & 4 & 3571.56 & 904 & 19 & 1623.54 & 554 & 18 & \textbf{1105.84} & 459 & 21 \\
50 & 5 & 4002.24 & 1023 & 25 & 1589.73 & 673 & 23 & \textbf{1179.27} & 483 & 19 \\ \hline
100 & 1 & $>$ 7200 & 1917 & 29 & 4731.48 & 1464 & 78 & \textbf{3430.76} & 1009 & 72 \\
100 & 2 & $>$ 7200 & 2532 & 21 & 3956.04 & 1576 & 36 & \textbf{3038.67} & 1162 & 26 \\
100 & 3 & $>$ 7200 & 2187 & 20 & 4121.46 & 1257 & 47 & \textbf{3267.31} & 1206 & 36 \\
100 & 4 & $>$ 7200 & 1375 & 43 & 4823.22 & 1042 & 56 & \textbf{3180.45} & 794 & 42 \\
100 & 5 & $>$ 7200 & 1967 & 19 & 3781.28 & 1272 & 31 & \textbf{2688.28} & 864 & 37 \\ \hline
150 & 1 & $>$ 7200 & 1688 & 32 & 6944.31 & 1110 & 94 & \textbf{5697.86} & 1043 & 83 \\
150 & 2 & $>$ 7200 & 1803 & 23 & 7019.38 & 1097 & 87 & \textbf{6246.27} & 976 & 79 \\
150 & 3 & $>$ 7200 & 1501 & 21 & 6759.04 & 1065 & 53 & \textbf{5393.86} & 979 & 40 \\
150 & 4 & $>$ 7200 & 1350 & 39 & 6491.46 & 993 & 55 & \textbf{5202.67} & 744 & 65 \\
150 & 5 & $>$ 7200 & 1635 & 33 & 6826.97 & 1046 & 81 & \textbf{5850.57} & 920 & 74 \\
\end{tabular}
}
\end{center}
\end{table}

\subsubsection{Comparison with Full Formulations} \label{subsubsec:full}
\blue{In this section, we compare the performance of the decomposition methods DD-BD and G-BD 2 with that of the full formulation \eqref{UCP:opj}--\eqref{const:continuous} when solved by Gurobi at its default settings.
Since it is well-known that decomposition methods are advantageous for larger problem sizes, we conduct these experiments for an increased number of demand scenarios, i.e., $|\Xi| \in \{200, 400, 600, 800\}$, to have a meaningful comparison.
These results are given in Tables~\ref{tab:2} and \ref{tab:3} for 20 and 50 generators, respectively. 
The solution time for larger number of generators exceeds the time limit of 7200 seconds for all methods.
In these tables, the first column shows the number of demand scenarios, and the second column contains the instance number.
The next three columns under ``G-BD 2" report the solution time, the number of feasibility and optimality cuts, respectively. 
Columns 6-8 contain similar results for the DD-BD method.
The last column shows the solution time for the full formulation.
}

\begin{figure}
\centering
\includegraphics[scale=0.60]{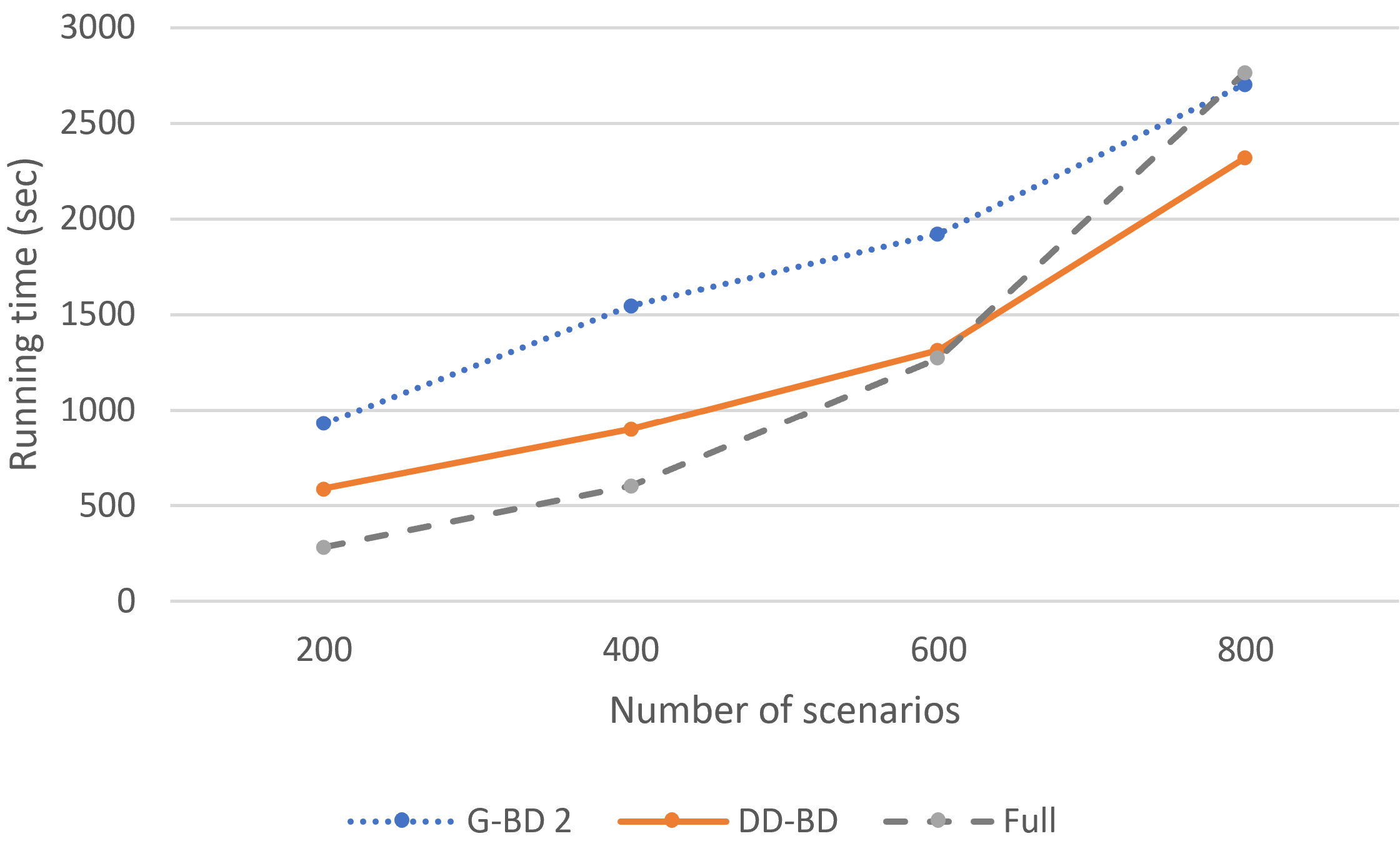}
\caption{\blue{Average solution times of the G-BD 2, DD-BD, and Full formulations for different number of scenarios and $n=20$}}
\label{fig:n20s}
\end{figure}

\begin{figure}
\centering
\includegraphics[scale=0.60]{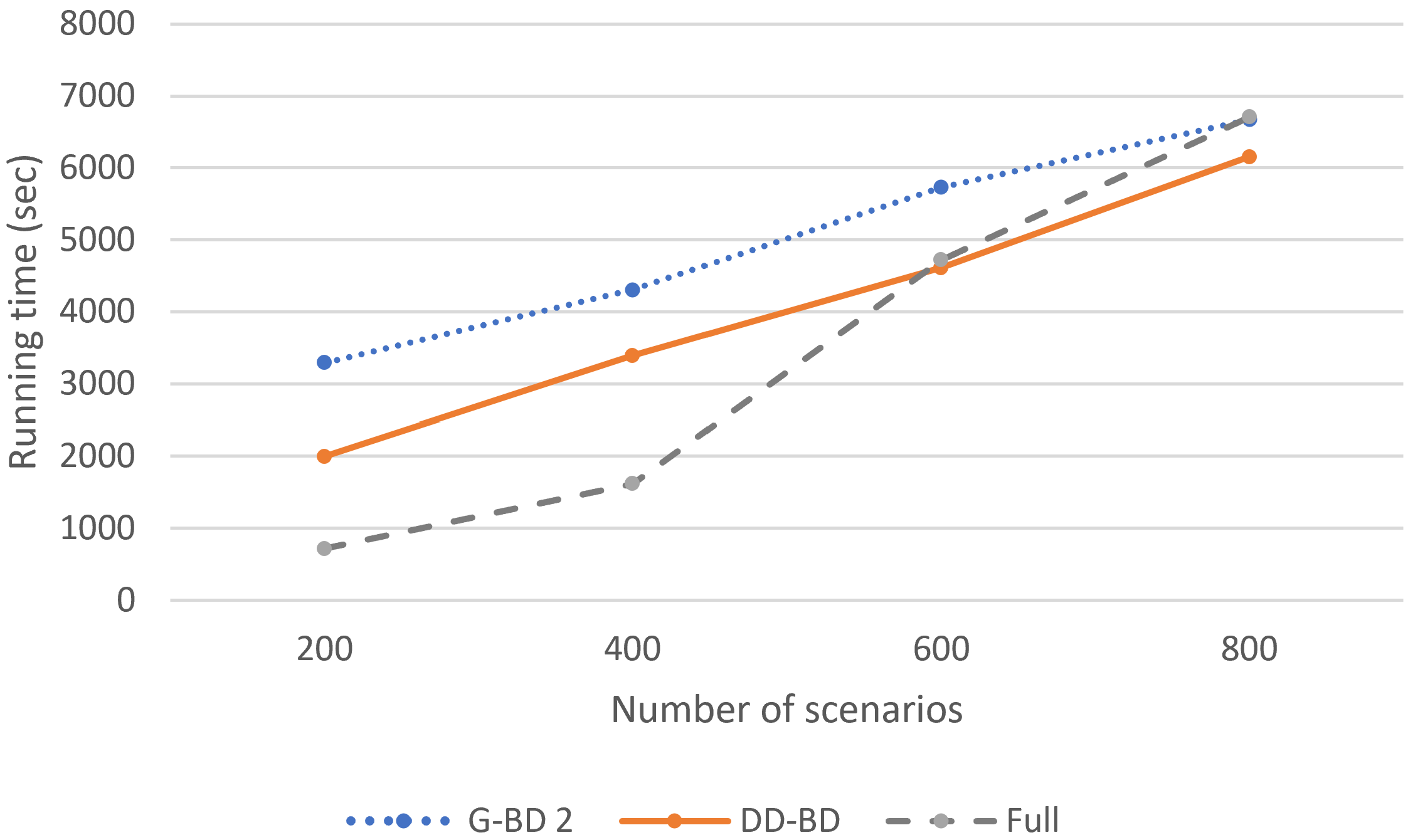}
\caption{\blue{Average solution times of the G-BD 2, DD-BD, and Full formulations for different number of scenarios and $n=50$}}
\label{fig:n50s}
\end{figure}

\blue{As evident in Tables~\ref{tab:2} and \ref{tab:3}, similarly to Table~\ref{tab:1}, the DD-BD method outperforms the G-BD 2 uniformly for all problem sizes.
Even though the full formulation solves the instances with 200 and 400 demand scenarios faster than the decomposition methods, it suffers from a steeper time increase as the number of demand scenarios increases. 
As a result, the DD-BD method catches up to the full formulation for 600 demand scenarios and overcomes it for the larger size with 800 demand scenarios, implying the superior performance of the DD-BD compared to other approaches for large problems sizes.
To better visualize this pattern, we present the average solution time for each size category against the number of demand scenarios in Figures~\ref{fig:n20s} and \ref{fig:n50s}.
}

\blue{We conclude this section by noting that the presented computational results are obtained through a basic implementation of the DD-BD without incorporating any computational enhancement methods for DDs such as \textit{variables ordering}, \textit{DD reduction}, \textit{dynamic width adjustment}, \textit{primal heuristics}, etc.
As common in the literature, such basic DD implementations are often compared against modern solvers with the presolve/cuts/heuristics settings turned off to provide a fair comparison for the core solution methodology.
Despite being at such disadvantage, the above results indicate that the proposed DD-BD method can still outperform the G-BD and the full UCP formulation solved with the solver's default settings.
While this shows the potential of the DD-BD method in energy applications, its performance can be further improved using modern computational practices in the DD literature, which can be considered as a direction for future work.}

\begin{table}[!htbp]
\begin{center}
\caption{\blue{Solution times (in seconds) of G-BD 2, DD-BD, and Full formulations for instances of the stochastic UCP with $n=20$}}
\label{tab:2}
\scalebox{1}{
\begin{tabular}{cc|rrr|rrr|rrr} 
\multicolumn{2}{c|}{} & \multicolumn{3}{c|}{{G-BD 2}} & \multicolumn{3}{c|}{{DD-BD}} & \multicolumn{1}{r}{{Full}}\\
$|\Xi|$ & Instance & time & f cuts & o cuts & time & f cuts & o cuts & time \\\hline
200 & 1 & 945.18 & 452 & 42 & 708.96 & 308 & 38 & \textbf{304.05} \\
200 & 2 & 1004.18 & 562 & 31 & 542.99 & 385 & 33 & \textbf{186.35} \\
200 & 3 & 1186.93 & 612 & 31 & 510.94 & 531 & 36 & \textbf{252.89} \\
200 & 4 & 809.82 & 550 & 42 & 537.76 & 379 & 44 & \textbf{290.98} \\
200 & 5 & 707.68 & 599 & 31 & 650.80 & 515 & 23 & \textbf{371.19} \\\hline
400 & 1 & 1618.80 & 1034 & 52 & 946.88 & 953 & 59 & \textbf{616.48} \\
400 & 2 & 1812.31 & 847 & 42 & 856.44 & 695 & 35 & \textbf{574.28} \\
400 & 3 & 1451.45 & 514 & 38 & 1025.23 & 395 & 32 & \textbf{598.08} \\
400 & 4 & 1158.77 & 564 & 58 & 942.89 & 600 & 41 & \textbf{606.82} \\
400 & 5 & 1696.87 & 712 & 54 & 721.38 & 575 & 61 & \textbf{618.77} \\\hline
600 & 1 & 2115.54 & 909 & 56 & 1357.43 & 736 & 54 & \textbf{1132.66} \\
600 & 2 & 1530.90 & 1163 & 53 & \textbf{1338.87} & 1314 & 48 & 1413.65 \\
600 & 3 & 1826.37 & 1241 & 60 & 1270.83 & 1024 & 48 & \textbf{1142.26} \\
600 & 4 & 2190.60 & 1226 & 80 & \textbf{1353.67} & 875 & 89 & 1488.27 \\
600 & 5 & 1950.40 & 984 & 68 & 1256.71 & 728 & 50 & \textbf{1192.87} \\\hline
800 & 1 & 2879.51 & 1410 & 76 & \textbf{2231.36} & 1049 & 61 & 2601.54 \\
800 & 2 & 2332.64 & 1119 & 70 & \textbf{2195.52} & 758 & 83 & 2653.63 \\
800 & 3 & 2651.88 & 1458 & 71 & \textbf{2254.60} & 1107 & 86 & 2921.57 \\
800 & 4 & 2958.79 & 1623 & 88 & \textbf{2288.74} & 1671 & 75 & 2772.84 \\
800 & 5 & 2707.54 & 1570 & 86 & \textbf{2635.46} & 1522 & 69 & 2897.46 \\
\end{tabular}
}
\end{center}
\end{table}

\begin{table}[!htbp]
\begin{center}
\caption{\blue{Solution times (in seconds) of G-BD 2, DD-BD, and Full formulations for instances of the stochastic UCP with $n=50$}}
\label{tab:3}
\scalebox{1}{
\begin{tabular}{cc|rrr|rrr|rrr} 
\multicolumn{2}{c|}{} & \multicolumn{3}{c|}{{G-BD 2}} & \multicolumn{3}{c|}{{DD-BD}} & \multicolumn{1}{r}{{Full}}\\
$|\Xi|$ & Instance & time & f cuts & o cuts & time & f cuts & o cuts & time \\\hline
200 & 1 & 3505.12 & 1643 & 78 & 1944.37 & 1549 & 90 & \textbf{610.70} \\
200 & 2 & 3797.28 & 1066 & 43 & 1816.58 & 970 & 47 & \textbf{712.24} \\
200 & 3 & 2824.65 & 972 & 30 & 2210.88 & 916 & 22 & \textbf{628.35} \\
200 & 4 & 2949.56 & 855 & 20 & 2034.54 & 664 & 15 & \textbf{850.44} \\
200 & 5 & 3375.60 & 1478 & 84 & 1981.47 & 1028 & 76 & \textbf{773.09} \\\hline
400 & 1 & 3859.64 & 1804 & 70 & 3324.22 & 1629 & 75 & \textbf{1752.74} \\
400 & 2 & 4526.40 & 860 & 72 & 3092.76 & 621 & 52 & \textbf{1939.73} \\
400 & 3 & 4129.86 & 1440 & 58 & 3920.40 & 1209 & 49 & \textbf{1658.43} \\
400 & 4 & 4653.48 & 1776 & 38 & 2971.42 & 1633 & 41 & \textbf{1463.31} \\
400 & 5 & 4367.37 & 1383 & 81 & 3656.57 & 1211 & 100 & \textbf{1277.31} \\\hline
600 & 1 & 5785.03 & 1251 & 109 & \textbf{4003.92} & 1075 & 91 & 4498.25 \\
600 & 2 & 5834.85 & 1537 & 99 & 4583.70 & 1618 & 87 & \textbf{4124.68} \\
600 & 3 & 5410.30 & 1401 & 60 & \textbf{4970.63} & 1162 & 51 & 5065.13 \\
600 & 4 & 5384.34 & 1950 & 53 & \textbf{4409.66} & 1482 & 62 & 4874.60 \\
600 & 5 & 5739.52 & 1136 & 89 & 5094.68 & 1056 & 67 & \textbf{5048.91} \\\hline
800 & 1 & 6396.84 & 2285 & 94 & \textbf{5963.94} & 1713 & 74 & 6883.02 \\
800 & 2 & 6575.52 & 2303 & 110 & \textbf{6204.16} & 1704 & 129 & 6840.31 \\
800 & 3 & 6606.54 & 1831 & 82 & \textbf{6353.15} & 1617 & 76 & 6643.29 \\
800 & 4 & 6861.92 & 1161 & 71 & \textbf{6065.85} & 940 & 59 & 6411.49 \\
800 & 5 & 6919.75 & 2843 & 107 & \textbf{6182.09} & 1998 & 121 & 6757.77
\end{tabular}
}
\end{center}
\end{table}

\section{Conclusion} \label{sec:conclusion}
In this paper, we propose a geometric decomposition that restructures a mixed integer set through a collection of hyper-rectangle formations. 
We show that this decomposition approach is the key to model optimization problems with both integer and continuous variables through decision diagrams. 
In this regard, we introduce a method, referred to as DD-BD, that extends the applicability domain of decision diagrams to MIPs. 
We evaluate the performance of DD-BD by applying it to the unit commitment problem in the electric grid market.

\section*{Acknowledgement} 
This project is sponsored in part by the Iowa Energy Center, Iowa Economic Development Authority and its utility partners.
We thank the anonymous referees and the Associate Editor for their helpful comments that contributed to improving the paper.


\bibliographystyle{informs2014} 
\bibliography{IPDD} 

\begin{thebibliography}{56}
\providecommand{\natexlab}[1]{#1}
\providecommand{\url}[1]{\texttt{#1}}
\providecommand{\urlprefix}{URL }

\bibitem[{Atakan et~al.(2017)Atakan, Lulli, \protect\BIBand{}
  Sen}]{atakan2017state}
Atakan S, Lulli G, Sen S (2017) A state transition mip formulation for the unit
  commitment problem. \emph{IEEE Transactions on Power Systems} 33(1):736--748.

\bibitem[{Baldick(1995)}]{baldick1995generalized}
Baldick R (1995) The generalized unit commitment problem. \emph{IEEE
  Transactions on Power Systems} 10(1):465--475.

\bibitem[{Bendotti et~al.(2018)Bendotti, Fouilhoux, \protect\BIBand{}
  Rottner}]{bendotti2018min}
Bendotti P, Fouilhoux P, Rottner C (2018) The min-up/min-down unit commitment
  polytope. \emph{Journal of Combinatorial Optimization} 36(3):1024--1058.

\bibitem[{Bendotti et~al.(2019)Bendotti, Fouilhoux, \protect\BIBand{}
  Rottner}]{bendotti2019complexity}
Bendotti P, Fouilhoux P, Rottner C (2019) On the complexity of the unit
  commitment problem. \emph{Annals of Operations Research} 274(1-2):119--130.

\bibitem[{Bergman \protect\BIBand{} Cire(2016)}]{bergman2016multiobjective}
Bergman D, Cire AA (2016) Multiobjective optimization by decision diagrams.
  \emph{International Conference on Principles and Practice of Constraint
  Programming}, 86--95 (Springer).

\bibitem[{Bergman \protect\BIBand{} Cire(2018)}]{bergman2018discrete}
Bergman D, Cire AA (2018) Discrete nonlinear optimization by state-space
  decompositions. \emph{Management Science} 64(10):4700--4720.

\bibitem[{Bergman et~al.(2015)Bergman, Cire, \protect\BIBand{} van
  Hoeve}]{bergman2015lagrangian}
Bergman D, Cire AA, van Hoeve WJ (2015) Lagrangian bounds from decision
  diagrams. \emph{Constraints} 20(3):346--361.

\bibitem[{Bergman et~al.(2016)Bergman, Cire, van Hoeve, \protect\BIBand{}
  Hooker}]{bergman:ci:va:ho:2016}
Bergman D, Cire AA, van Hoeve WJ, Hooker J (2016) \emph{Decision Diagrams for
  Optimization} (Springer International Publishing).

\bibitem[{Bertsimas et~al.(2012)Bertsimas, Litvinov, Sun, Zhao,
  \protect\BIBand{} Zheng}]{bertsimas2012adaptive}
Bertsimas D, Litvinov E, Sun XA, Zhao J, Zheng T (2012) Adaptive robust
  optimization for the security constrained unit commitment problem. \emph{IEEE
  transactions on power systems} 28(1):52--63.

\bibitem[{Blanco \protect\BIBand{} Morales(2017)}]{blanco2017efficient}
Blanco I, Morales JM (2017) An efficient robust solution to the two-stage
  stochastic unit commitment problem. \emph{IEEE Transactions on Power Systems}
  32(6):4477--4488.

\bibitem[{Bonami et~al.(2020)Bonami, Salvagnin, \protect\BIBand{}
  Tramontani}]{bonami2020implementing}
Bonami P, Salvagnin D, Tramontani A (2020) Implementing automatic benders
  decomposition in a modern mip solver. \emph{International conference on
  integer programming and combinatorial optimization}, 78--90 (Springer).

\bibitem[{Carri{\'o}n \protect\BIBand{}
  Arroyo(2006)}]{carrion2006computationally}
Carri{\'o}n M, Arroyo JM (2006) A computationally efficient mixed-integer
  linear formulation for the thermal unit commitment problem. \emph{IEEE
  Transactions on power systems} 21(3):1371--1378.

\bibitem[{Codato \protect\BIBand{} Fischetti(2006)}]{codato2006combinatorial}
Codato G, Fischetti M (2006) Combinatorial benders' cuts for mixed-integer
  linear programming. \emph{Operations Research} 54(4):756--766.

\bibitem[{Conforti et~al.(2014{\natexlab{a}})Conforti, Cornu{\' e}jols,
  \protect\BIBand{} Zambelli}]{conforti:co:za:2014}
Conforti M, Cornu{\' e}jols G, Zambelli G (2014{\natexlab{a}}) \emph{Integer
  Programming} (Springer).

\bibitem[{Conforti et~al.(2014{\natexlab{b}})Conforti, Cornu{\'e}jols, Zambelli
  et~al.}]{conforti2014integer}
Conforti M, Cornu{\'e}jols G, Zambelli G, et~al. (2014{\natexlab{b}})
  \emph{Integer programming}, volume 271 (Springer).

\bibitem[{Davarnia(2021)}]{davarnia2020strong}
Davarnia D (2021) Strong relaxations for continuous nonlinear programs based on
  decision diagrams. \emph{Operations Research Letters}
  \urlprefix\url{http://dx.doi.org/10.1016/j.orl.2021.01.011}.

\bibitem[{Davarnia \protect\BIBand{} van Hoeve(2020)}]{davarnia:va:2020}
Davarnia D, van Hoeve WJ (2020) Outer approximation for integer nonlinear
  programs via decision diagrams. \emph{Mathematical Programming}
  \urlprefix\url{http://dx.doi.org/10.1007/s10107-020-01475-4}.

\bibitem[{Feng \protect\BIBand{} Ryan(2016)}]{feng2016solution}
Feng Y, Ryan SM (2016) Solution sensitivity-based scenario reduction for
  stochastic unit commitment. \emph{Computational Management Science}
  13(1):29--62.

\bibitem[{Fischetti et~al.(2010)Fischetti, Salvagnin, \protect\BIBand{}
  Zanette}]{fischetti2010note}
Fischetti M, Salvagnin D, Zanette A (2010) A note on the selection of
  benders’ cuts. \emph{Mathematical Programming} 124(1):175--182.

\bibitem[{Frangioni \protect\BIBand{} Gentile(2006)}]{frangioni2006solving}
Frangioni A, Gentile C (2006) Solving nonlinear single-unit commitment problems
  with ramping constraints. \emph{Operations Research} 54(4):767--775.

\bibitem[{Frangioni et~al.(2008)Frangioni, Gentile, \protect\BIBand{}
  Lacalandra}]{frangioni2008tighter}
Frangioni A, Gentile C, Lacalandra F (2008) Tighter approximated milp
  formulations for unit commitment problems. \emph{IEEE Transactions on Power
  Systems} 24(1):105--113.

\bibitem[{Franz et~al.(2020)Franz, Rieck, \protect\BIBand{}
  Zimmermann}]{franz2020long}
Franz A, Rieck J, Zimmermann J (2020) A long-term unit commitment problem with
  hydrothermal coordination for economic and emission control in large-scale
  electricity systems. \emph{OR spectrum} 42(1):235--259.

\bibitem[{Fu et~al.(2013)Fu, Li, \protect\BIBand{} Wu}]{fu2013modeling}
Fu Y, Li Z, Wu L (2013) Modeling and solution of the large-scale
  security-constrained unit commitment. \emph{IEEE Transactions on Power
  Systems} 28(4):3524--3533.

\bibitem[{Garver(1962)}]{garver1962power}
Garver LL (1962) Power generation scheduling by integer
  programming-{D}evelopment of theory. \emph{Transactions of the American
  Institute of Electrical Engineers. Part III: Power Apparatus and Systems}
  81(3):730--734.

\bibitem[{Gonzalez et~al.(2020)Gonzalez, Cire, Lodi, \protect\BIBand{}
  Rousseau}]{gonzalez:ci:lo:ro:2020}
Gonzalez J, Cire A, Lodi A, Rousseau LM (2020) Integrated integer programming
  and decision diagram search tree with an application to the maximum
  independent set problem. \emph{Constraints} 25:23--46.

\bibitem[{Guan et~al.(2003)Guan, Zhai, \protect\BIBand{}
  Papalexopoulos}]{guan2003optimization}
Guan X, Zhai Q, Papalexopoulos A (2003) Optimization based methods for unit
  commitment: Lagrangian relaxation versus general mixed integer programming.
  \emph{2003 IEEE Power Engineering Society General Meeting (IEEE Cat. No.
  03CH37491)}, volume~2, 1095--1100 (IEEE).

\bibitem[{Had{\u z}i{\' c} \protect\BIBand{} Hooker(2006)}]{hadzic:ho:2006}
Had{\u z}i{\' c} T, Hooker JN (2006) Discrete global optimization with binary
  decision diagrams. \emph{Workshop on Global Optimization: Integrating
  Convexity, Optimization, Logic Programming, and Computational Algebraic
  Geormetry (GICOLAG)}.

\bibitem[{Huang et~al.(2017)Huang, Pardalos, \protect\BIBand{}
  Zheng}]{huang2017electrical}
Huang Y, Pardalos PM, Zheng QP (2017) \emph{Electrical power unit commitment:
  deterministic and two-stage stochastic programming models and algorithms}
  (Springer).

\bibitem[{Knueven et~al.(2020{\natexlab{a}})Knueven, Ostrowski,
  \protect\BIBand{} Watson}]{knueven2020novel}
Knueven B, Ostrowski J, Watson JP (2020{\natexlab{a}}) A novel matching
  formulation for startup costs in unit commitment. \emph{Mathematical
  Programming Computation} 1--24.

\bibitem[{Knueven et~al.(2020{\natexlab{b}})Knueven, Ostrowski,
  \protect\BIBand{} Watson}]{knueven2020mixed}
Knueven B, Ostrowski J, Watson JP (2020{\natexlab{b}}) On mixed-integer
  programming formulations for the unit commitment problem. \emph{INFORMS
  Journal on Computing} .

\bibitem[{Lee \protect\BIBand{} Feng(1992)}]{lee1992multi}
Lee FN, Feng Q (1992) Multi-area unit commitment. \emph{IEEE Transactions on
  power systems} 7(2):591--599.

\bibitem[{Li \protect\BIBand{} Shahidehpour(2005)}]{li2005price}
Li T, Shahidehpour M (2005) Price-based unit commitment: A case of lagrangian
  relaxation versus mixed integer programming. \emph{IEEE transactions on power
  systems} 20(4):2015--2025.

\bibitem[{Li et~al.(2007)Li, McCalley, \protect\BIBand{} Ryan}]{li2007risk}
Li Y, McCalley JD, Ryan S (2007) Risk-based unit commitment. \emph{2007 IEEE
  Power Engineering Society General Meeting}, 1--7 (IEEE).

\bibitem[{Lorca et~al.(2016)Lorca, Sun, Litvinov, \protect\BIBand{}
  Zheng}]{lorca2016multistage}
Lorca {\'A}, Sun XA, Litvinov E, Zheng T (2016) Multistage adaptive robust
  optimization for the unit commitment problem. \emph{Operations Research}
  64(1):32--51.

\bibitem[{Maifeld \protect\BIBand{} Sheble(1996)}]{maifeld1996genetic}
Maifeld TT, Sheble GB (1996) Genetic-based unit commitment algorithm.
  \emph{IEEE Transactions on Power systems} 11(3):1359--1370.

\bibitem[{Mantawy et~al.(1998)Mantawy, Abdel-Magid, \protect\BIBand{}
  Selim}]{mantawy1998simulated}
Mantawy A, Abdel-Magid YL, Selim SZ (1998) A simulated annealing algorithm for
  unit commitment. \emph{IEEE Transactions on Power Systems} 13(1):197--204.

\bibitem[{Morales-Espa{\~n}a et~al.(2013)Morales-Espa{\~n}a, Latorre,
  \protect\BIBand{} Ramos}]{morales2013tight}
Morales-Espa{\~n}a G, Latorre JM, Ramos A (2013) Tight and compact milp
  formulation for the thermal unit commitment problem. \emph{IEEE Transactions
  on Power Systems} 28(4):4897--4908.

\bibitem[{Nemhauser \protect\BIBand{} Wolsey(1990)}]{nemhauser1990recursive}
Nemhauser GL, Wolsey LA (1990) A recursive procedure to generate all cuts for
  0--1 mixed integer programs. \emph{Mathematical Programming}
  46(1-3):379--390.

\bibitem[{Ostrowski et~al.(2011)Ostrowski, Anjos, \protect\BIBand{}
  Vannelli}]{ostrowski2011tight}
Ostrowski J, Anjos MF, Vannelli A (2011) Tight mixed integer linear programming
  formulations for the unit commitment problem. \emph{IEEE Transactions on
  Power Systems} 27(1):39--46.

\bibitem[{Padhy(2004)}]{padhy2004unit}
Padhy NP (2004) Unit commitment-a bibliographical survey. \emph{IEEE
  Transactions on power systems} 19(2):1196--1205.

\bibitem[{Pang et~al.(1981)Pang, Shebl{\'e}, \protect\BIBand{}
  Albuyeh}]{pang1981evaluation}
Pang C, Shebl{\'e} GB, Albuyeh F (1981) Evaluation of dynamic programming based
  methods and multiple area representation for thermal unit commitments.
  \emph{IEEE Transactions on Power Apparatus and Systems} (3):1212--1218.

\bibitem[{Queyranne \protect\BIBand{} Wolsey(2017)}]{queyranne2017tight}
Queyranne M, Wolsey LA (2017) Tight mip formulations for bounded up/down times
  and interval-dependent start-ups. \emph{Mathematical Programming}
  164(1-2):129--155.

\bibitem[{Rajan et~al.(2003)Rajan, Mohan, \protect\BIBand{}
  Manivannan}]{rajan2003neural}
Rajan CA, Mohan M, Manivannan K (2003) Neural-based tabu search method for
  solving unit commitment problem. \emph{IEE Proceedings-Generation,
  Transmission and Distribution} 150(4):469--474.

\bibitem[{Rajan et~al.(2005)Rajan, Takriti et~al.}]{rajan2005minimum}
Rajan D, Takriti S, et~al. (2005) Minimum up/down polytopes of the unit
  commitment problem with start-up costs. \emph{IBM Res. Rep} 23628:1--14.

\bibitem[{Rockafeller(1970)}]{rockafellar:1970}
Rockafeller RT (1970) \emph{Convex Analysis} (New Jersey, NJ: Princeton
  University Press).

\bibitem[{Rodr{\'\i}guez et~al.(2021)Rodr{\'\i}guez, Anjos, C{\^o}t{\'e},
  \protect\BIBand{} Desaulniers}]{rodriguez2021accelerating}
Rodr{\'\i}guez JA, Anjos MF, C{\^o}t{\'e} P, Desaulniers G (2021) Accelerating
  benders decomposition for short-term hydropower maintenance scheduling.
  \emph{European Journal of Operational Research} 289(1):240--253.

\bibitem[{Saneifard et~al.(1997)Saneifard, Prasad, \protect\BIBand{}
  Smolleck}]{saneifard1997fuzzy}
Saneifard S, Prasad NR, Smolleck HA (1997) A fuzzy logic approach to unit
  commitment. \emph{IEEE Transactions on Power Systems} 12(2):988--995.

\bibitem[{Serra \protect\BIBand{} Hooker(2020)}]{serra:ho:2020}
Serra T, Hooker J (2020) Compact representation of near-optimal integer
  programming solutions. \emph{Mathematical Programming} 182:199--232.

\bibitem[{Silbernagl et~al.(2015)Silbernagl, Huber, \protect\BIBand{}
  Brandenberg}]{silbernagl2015improving}
Silbernagl M, Huber M, Brandenberg R (2015) Improving accuracy and efficiency
  of start-up cost formulations in mip unit commitment by modeling power plant
  temperatures. \emph{IEEE Transactions on Power Systems} 31(4):2578--2586.

\bibitem[{Tuffaha \protect\BIBand{} Gravdahl(2013)}]{tuffaha2013mixed}
Tuffaha M, Gravdahl JT (2013) Mixed-integer formulation of unit commitment
  problem for power systems: Focus on start-up cost. \emph{IECON 2013-39th
  Annual Conference of the IEEE Industrial Electronics Society}, 8160--8165
  (IEEE).

\bibitem[{van Ackooij et~al.(2018)van Ackooij, Danti~Lopez, Frangioni,
  Lacalandra, \protect\BIBand{} Tahanan}]{van2018large}
van Ackooij W, Danti~Lopez I, Frangioni A, Lacalandra F, Tahanan M (2018)
  Large-scale unit commitment under uncertainty: an updated literature survey.
  \emph{Annals of Operations Research} 271(1):11--85.

\bibitem[{van~der Linden(2017)}]{van2017decision}
van~der Linden K (2017) Decision diagrams for decomposed mixed integer linear
  programs. Master's thesis. Delft University of Technology.

\bibitem[{Wu(2011)}]{wu2011tighter}
Wu L (2011) A tighter piecewise linear approximation of quadratic cost curves
  for unit commitment problems. \emph{IEEE Transactions on Power Systems}
  26(4):2581--2583.

\bibitem[{Xavier et~al.(2020)Xavier, Kazachkov, \protect\BIBand{}
  Qiu}]{UnitCommitmentjl}
Xavier AS, Kazachkov AM, Qiu F (2020) Unitcommitment.jl: A julia/jump
  optimization package for security-constrained unit commitment. Zenodo (2020).
  DOI: 10.5281/zenodo.4269874.

\bibitem[{Zheng et~al.(2014)Zheng, Wang, \protect\BIBand{}
  Liu}]{zheng2014stochastic}
Zheng QP, Wang J, Liu AL (2014) Stochastic optimization for unit commitment—a
  review. \emph{IEEE Transactions on Power Systems} 30(4):1913--1924.

\bibitem[{Zheng et~al.(2013)Zheng, Wang, Pardalos, \protect\BIBand{}
  Guan}]{zheng2013decomposition}
Zheng QP, Wang J, Pardalos PM, Guan Y (2013) A decomposition approach to the
  two-stage stochastic unit commitment problem. \emph{Annals of Operations
  Research} 210(1):387--410.

\end{thebibliography}

\clearpage

\begin{APPENDICES}
\section{Omitted Proofs} \label{sec:proofs}
\proof{\textbf{Proof of Theorem~\ref{thm:equivalence_general}}.}
\blue{For a given function $f(\vc{x})$ that is convex in $\vc{x}_I$}, define $z^1 = \max\{f(\vc{x}) | \vc{x} \in \bigcup_{j \in J} \mc{X}(P_I^j)\}$ and $z^2 = \max\{f(\vc{x}) | \vc{x} \in \bigcup_{j \in J} \bigcup_{k \in K_j} R_j^k\}$.
We first show that $z^1 = z^2$.
For each $j \in J$, we write that $\mc{X}(P_I^j) \subseteq \mc{X}(\bigcup_{k \in K_j} R_j^k) \subseteq \bigcup_{k \in K_j} R_j^k$, where the first inclusion follows from condition (iii), and the second inclusion follows from definition of extreme points.
Therefore, $z^1 \leq z^2$.
The above inclusions also show that $\bigcup_{j \in J} \mc{X}(P_I^j)$ is a finite set since $\bigcup_{j \in J} \mc{X}(\bigcup_{k \in K_j} R_j^k)$ is finite because of the finiteness of $J$, $K_j$ and the set of extreme points of hyper-rectangles $R_j^k$.
Let $\vc{x}^*$ be an optimal solution for $z^2$.
\blue{Such optimal solution exists because of the compactness of the feasible region.}
It follows that $\vc{x}^* \in \bigcup_{k \in K_{j^*}} R_{j^*}^k$ for some $j^* \in J$.
\blue{Condition (i) implies that coordinates $i \in N \setminus I$ are fixed in $P_I^{j^*}$, and condition (iii) implies that these coordinates must also be fixed for $\bigcup_{k \in K_{j^*}} R_{j^*}^k$.}
Therefore, $x^*_i$ is fixed at coordinates $i \in N \setminus I$.
Since $f(\vc{x})$ is convex in the unfixed variables $\vc{x}_I$, its maximum over $\bigcup_{k \in K_{j^*}} R_{j^*}^k$ occurs at an extreme point $\vc{\bar{x}}$ of $\conv(\bigcup_{k \in K_{j^*}}R_{j^*}^k)$. Condition (iii) implies that $\vc{\bar{x}} \in \mc{X}(P_I^{j^*})$, proving that $z^1 \geq z^2$.
For the second part of the proof, we show that $z^3 = \max\{f(\vc{x}) | \vc{x} \in \mc{P}\} = z^1 = z^2$.
It follows from condition (ii) that $z^1 \leq z^3 \leq z^4$ where $z^4 = \max\{f(\vc{x}) | \vc{x} \in \bigcup_{j \in J} P_I^j\}$.
Using an argument similar to that given above, we conclude that the optimal value $z^4$ is attained at an extreme point of $P_I^{j^*}$ for some $j^* \in J$, which is also an extreme point of $\bigcup_{k \in K_{j^*}}R_{j^*}^k$.
Therefore, \blue{using the definition of $z^4$ and $z^2$}, we write that $z^4 \leq z^2 = z^1$, proving the result.
\Halmos
\endproof
\medskip
\proof{\textbf{Proof of Proposition~\ref{prop:converse}}.}
    Consider any point $\bar{\vc{x}}_{N \setminus I} \in \proj_{x_{N \setminus I}}(\mc{P})$, and define the indicator function 
    \[\delta(\vc{x}|\bar{\vc{x}}_{N \setminus I}) = \left\{ \begin{array}{ll}
        0, & \text{ if } \vc{x}_{N \setminus I} = \bar{\vc{x}}_{N \setminus I} \\
        -\infty, & \text{ else. }
    \end{array}  \right.\]    
    We can use $\delta(\vc{x}|\bar{\vc{x}}_{N \setminus I})$ as the objective function in the relation $\max\{f(\vc{x}) \suchthat \vc{x} \in \mc{P}\} = \max\{f(\vc{x}) \suchthat \vc{x} \in \mc{Q}\}$ as it is convex in $\vc{x}_I$.
    \blue{This follows from the fact that $\delta(\vc{x}|\bar{\vc{x}}_{N \setminus I}) = 0$ when $\vc{x}_{N \setminus I} = \bar{\vc{x}}_{N \setminus I}$, and it is $-\infty$ otherwise, which is \textit{improper} convex; see \cite{rockafellar:1970}.} 
    We conclude that $\bar{\vc{x}}_{N \setminus I} \in \proj_{x_{N \setminus I}}(\mc{Q})$.
    Using a similar argument for the reverse direction, we conclude that $\proj_{x_{N \setminus I}}(\mc{P}) = \proj_{x_{N \setminus I}}(\mc{Q})$.
    This also shows that $\proj_{x_{N \setminus I}}(\mc{P})$ is finite as $\mc{Q}$ is finite \blue{by assumption}.
    For any point $\bar{\vc{x}}_{N \setminus I} \in \proj_{x_{N \setminus I}}(\mc{P})$, define $\mc{P}(\bar{\vc{x}}_{N \setminus I}) = \mc{P}\cap \{\vc{x} | \vc{x}_{N \setminus I} = \bar{\vc{x}}_{N \setminus I}\}$, and $\mc{Q}(\bar{\vc{x}}_{N \setminus I}) = \mc{Q}\cap \{\vc{x} | \vc{x}_{N \setminus I} = \bar{\vc{x}}_{N \setminus I}\}$.
    For any convex function $f_I(\vc{x}_I)$ defined in the space of variables $\vc{x}_I$, and any point $\bar{\vc{x}}_{N \setminus I} \in \proj_{x_{N \setminus I}}(\mc{P})$, we have that 
    \begin{subequations}
    \begin{align}
    	& \max\{f_I(\vc{x}_I) | \vc{x} \in \mc{P}(\bar{\vc{x}}_{N \setminus I})\} \label{eq:chain1} \\
    	= & \max\{\delta(\vc{x}|\bar{\vc{x}}_{N \setminus I}) + f_I(\vc{x}_I) | \vc{x} \in \mc{P}\} \label{eq:chain2} \\
    	= & \max\{\delta(\vc{x}|\bar{\vc{x}}_{N \setminus I}) + f_I(\vc{x}_I) | \vc{x} \in \mc{Q}\} \label{eq:chain3} \\
    	= & \max\{f_I(\vc{x}_I) | \vc{x} \in \mc{Q}(\bar{\vc{x}}_{N \setminus I})\}, \label{eq:chain4}
    \end{align}
    \end{subequations}
    where \eqref{eq:chain1} and \eqref{eq:chain4} follow from the definition of the indicator function above, and \eqref{eq:chain2} follows from the assumption.
    Since $f_I(\vc{x}_I)$ is convex, the optimal value of \eqref{eq:chain1} is attained at an extreme point $\hat{\vc{x}}$ of $\mc{P}(\bar{\vc{x}}_{N \setminus I})$.
    For each such extreme point, there are infinitely many convex functions $\hat{f}_I(\vc{x}_I)$ with $\hat{\vc{x}}$ as their unique maximizer over $\mc{P}(\bar{\vc{x}}_{N \setminus I})$.
    For every one of these functions, according to the chain equalities \eqref{eq:chain1}--\eqref{eq:chain4}, $\mc{Q}$ has a point that matches the optimal value.
    Since $\mc{Q}$ is finite, it must be that $\hat{\vc{x}} \in \mc{Q}(\bar{\vc{x}}_{N \setminus I})$.
    This shows that the set of extreme points of $\mc{P}(\bar{\vc{x}}_{N \setminus I})$ is finite.
    Similarly, it can be shown that for any extreme point $\dot{\vc{x}}$ of $\conv(\mc{Q}(\bar{\vc{x}}_{N \setminus I}))$, we have that $\dot{\vc{x}} \in \mc{P}(\bar{\vc{x}}_{N \setminus I})$.
    As a result, $\conv(\mc{P}(\bar{\vc{x}}_{N \setminus I})) = \conv(\mc{Q}(\bar{\vc{x}}_{N \setminus I}))$ for every $\bar{\vc{x}}_{N \setminus I} \in \proj_{x_{N \setminus I}}(\mc{P})$.
    Now, construct sets $P_I^j$ in Theorem~\ref{thm:equivalence_general} as $\conv(\mc{P}(\bar{\vc{x}}_{N \setminus I}))$ for every $\bar{\vc{x}}_{N \setminus I} \in \proj_{x_{N \setminus I}}(\mc{P})$, which yields a finite collection.
    Condition (i) is satisfied as each set $P_I^j$ is restricted at $\{\vc{x} | \vc{x}_{N \setminus I} = \bar{\vc{x}}_{N \setminus I}\}$.
    Condition (ii) holds since for any extreme point of $P_I^j$, there is a point of $\mc{P}$ by construction, and since any point $\bar{\vc{x}} \in \mc{P}$ satisfies $\bar{\vc{x}} \in \mc{P}(\bar{\vc{x}}_{N \setminus I}) \subseteq \conv(\mc{P}(\bar{\vc{x}}_{N \setminus I}))$.
    For condition (iii), rectangles $R_j^k$ can be considered as the set of extreme points of $P_I^j$, which has been shown to be finite. 
 \Halmos
\endproof
\medskip
\proof{\textbf{Proof of Lemma~\ref{lem:equivalence}}.}
The result is obtained as a special case of Theorem~\ref{thm:equivalence_general}, where sets $P_I^j$ and $R_j^k$ coincide, i.e., $P_I^j = R_j^1$ for all $j \in J$. 
The conditions and the result follow immediately.
\Halmos
\endproof
\medskip
\proof{\textbf{Proof of Proposition~\ref{prop:reduction}}.}
For \blue{a given DD} $\mc{D}$, define a node-sequence as an ordered set of connected nodes from the root to the terminal, i.e., $\vc{u} = (u_1, u_2, \cdots, u_{n+1})$ where $u_i \in \mc{U}_i$ for $i \in N \cup \{n+1\}$.
Let $U$ be the collection of all node-sequences of $\mc{D}$.
For $\vc{u} \in U$, define 
\begin{equation*}
S_I(\vc{u}) = \Bigg\{x \in \Re^n \, \Bigg| \, \begin{array}{ll}
l_{(u_i,u_{i+1})}^{\text{min}} \leq x_i \leq l_{(u_i,u_{i+1})}^{\text{max}}, & \forall i \in I \\
x_i = l_{(u_i,u_{i+1})}, & \forall i \in N \setminus I 
\end{array}  \Bigg\}.
\end{equation*}
Viewing $\Sol(\mc{D})$ as \blue{a compact }set $\mc{P}$ in Lemma~\ref{lem:equivalence}, it is straightforward to verify that sets $S_I(\vc{u})$ satisfy the conditions for hyper-rectangles $R_I^j$.
It also follows from the definition of (virtual) DDs that $\Sol(\hat{\mc{D}}) = \bigcup_{\vc{u} \in U} S_I(\vc{u})$ and $\Sol(\bar{\mc{D}}) = \bigcup_{\vc{u} \in U} \mc{X}(S_I(\vc{u}))$.
The result follows from Lemma~\ref{lem:equivalence}.  
\Halmos
\endproof
\medskip
\proof{\textbf{Proof of Corollary~\ref{cor:cont_DD}}.}
    For the direct implication, assume that \blue{a given compact set} $\mc{P}$ admits a rectangular decomposition w.r.t. $I$ through sets $P_I^j$ for $j \in J$.
    Then, the DD that encodes the finite collection of points $\bigcup_{j \in J} \mc{X}(P_I^j)$ provides the desired DD representation because of Theorem~\ref{thm:equivalence_general}.
    For the reverse implication, assume that \blue{a given compact set }$\mc{P}$ is DD-representable w.r.t. $I$ through \blue{a DD }$\mc{D}$.
    Since $\Sol(\mc{D})$ is finite, it follows from Proposition~\ref{prop:converse} that $\mc{P}$ admits a rectangular decomposition w.r.t. $I$.
	\Halmos
\endproof
\medskip
\proof{\textbf{Proof of Corollary~\ref{cor:mip}}.}
	Since $\mc{Q}$ is bounded, there are finitely many points $\bar{\vc{x}} \in \proj_x(\mc{P})$.
	For each such point, it follows from the definition of projection and the boundedness of $\mc{Q}$ that there exists an interval $[l_{\bar{x}}, u_{\bar{x}}]$ such that $(\bar{\vc{x}}; \bar{y}) \in \mc{P}$ for every $\bar{y} \in [l_{\bar{x}}, u_{\bar{x}}]$.
	As a result, we can write that $\mc{P} = \bigcup_{\bar{\vc{x}} \in \proj_x(\mc{P})} R_y^{\bar{x}}$ where $R_y^{\bar{x}} = \{(\vc{x};y) \in \Re^{n+1} \, | \, \vc{x} = \bar{\vc{x}}, \, y \in [l_{\bar{x}}, u_{\bar{x}}]\}$.
	Hyper-rectangles $R_y^{\bar{x}}$ satisfy the conditions of Lemma~\ref{lem:equivalence}, and hence admits a rectangular decomposition w.r.t. the index of variable $y$ which is $n+1$.
	The result follows from Corollary~\ref{cor:cont_DD}.
	\Halmos
\endproof
\medskip
\proof{\textbf{Proof of Theorem~\ref{thm:DD-BD}}.} 
First, we show that the Algorithm~\ref{alg:DD-BD} terminates after a finite number of iterations.
To this end, we argue that each loop in the algorithm is repeated for a finite number of iterations.
The outer \textit{while} loop is executed for each member of the partial assignment set $\hat{\mc{X}}$.
These partial assignments are generated based on longest $r$-$u$ paths associated with nodes $u$ in the exact cut set of relaxed DDs $\overline{\mc{D}}$, which contain a finite number of paths by definition; see line 22 of the algorithm.
Each resulting partial assignment $\hat{\vc{x}}$ is different due to the structure of DDs that do not admit multiple paths with similar encoding values.
As a result, an upper bound for the number of partial assignments that can be included in $\hat{\mc{X}}$ is all possible partial assignments of the solution set of $\mc{M}$ for $\vc{x}$ variables.
Because the discrete set $\mc{P}$ in the description of $\mc{M}$ is bounded by assumption, we conclude that $\hat{\mc{X}}$ contains a finite number of elements. 
For the \textit{repeat} loop in lines 6--10, the goal is to obtain the optimal value of the solution set represented by $\underline{\mc{D}}$ subject to the constraints of the subproblem $\mc{S}(.)$.
It follows from the property of BD method that the resulting optimality and feasibility cuts correspond to the extreme points and extreme rays of the dual of the subproblem, which is independent of the choice of the fixed point $\underline{\vc{x}}$.
Since the feasible region of this dual problem is a polyhedron, the set of its extreme points and rays are finite, which yield a finite set of optimality and feasibility cuts that can be added through the execution of this loop.
\blue{It is also known in the BD literature (\cite{conforti:co:za:2014}) that the same optimality/feasibility cuts cannot be generated repeatedly through different iterations, since the added cuts remain in the description of the master DD due to refinement.}
A similar argument can be used to show that the number of iterations of the \textit{repeat} loop in lines 18--21 is finite, thereby yielding the result.

Next, we show that the outputs $(\vc{x}^*, z^*)$ and $w^*$ of Algorithm~\ref{alg:DD-BD} give an optimal solution and optimal value of $\mc{H}$, respectively.
To this end, we first prove that $(\vc{x}^*, z^*)$ is a feasible solution to $\mc{H}$ and $w^*$ is a lower bound for its optimal value.
This solution is updated at line 12 of Algorithm~\ref{alg:DD-BD} as the point encoding a longest $r$-$t$ path of $\underline{\mc{D}}$, after being refined with respect to optimality and feasibility cuts generated from subproblems $\mc{S}(.)$. 
We write that $(\vc{x}^*, z^*) \in \Sol(\underline{\mc{D}}) \subseteq \mc{M}^C(\hat{\vc{x}}) \subseteq \mc{M}$, where the second inclusion holds because $\underline{\mc{D}}$ is a restricted DD associated with $\mc{M}^C(\hat{\vc{x}})$ for some $C$ and $\hat{\vc{x}}$, and the third inclusion follows from the fact that a feasible solution to $\mc{M}^C(\hat{\vc{x}})$ is feasible to $\mc{M}$ by definition.
Further, the BD \blue{structure} implies that the point encoding a longest path obtained at the termination of the loop in line 6--10 satisfies the constraints of the subproblem $\mc{S}(.)$.
As a result, $(\vc{x}^*, z^*)$ satisfies the constraints of both the master and the subproblem, hence being feasible to $\mc{H}$.
Using a similar argument, we obtain that $w^*$ is a lower bound for the optimal value of $\mc{M}$ subject to constraints in $\mc{S}(.)$, since $w^*$ is the length of the longest path in the restricted DD representing a restriction of $\mc{M}$, after valid optimality and feasibility cuts are added based on the subproblem.
Now, we show that $(\vc{x}^*, z^*)$ is indeed an optimal solution to $\mc{H}$.
Assume by contradiction that there exists an optimal solution $(\tilde{\vc{x}}, \tilde{z})$ of $\mc{H}$ with optimal value $\tilde{w}$ such that $\tilde{w} > w^*$.
There are three cases.
(i) Assume that $\tilde{\vc{x}}$ is added as a partial assignment to set $\hat{\mc{X}}$ at some iteration of the algorithm.
For the while loop where $\tilde{\vc{x}}$ is selected at line 3, all $\vc{x}$ variables are fixed for the restricted DD $\underline{\mc{D}}$.
Therefore, $\tilde{\vc{x}}$ is used as the input for the subproblem, i.e., we solve $\mc{S}(\tilde{\vc{x}})$, which yields the optimal value $\tilde{z}$ since $(\tilde{\vc{x}}, \tilde{z})$ is an optimal solution of $\mc{H}$.
Since the weights on the arcs of the restricted DD $\underline{\mc{D}}$ are set as the coefficients of variables in the linear objective function of $\mc{H}$, the length of the $r$-$t$ path of $\underline{\mc{D}}$ encoded by $(\tilde{\vc{x}}, \tilde{z})$ is $\tilde{w}$.
If follows from the contradiction assumption that $\tilde{w} > w^*$ in line 11 of Algorithm~\ref{alg:DD-BD}, which leads to updating the optimal solution and optimal value to $(\tilde{\vc{x}}, \tilde{z})$ and $\tilde{w}$.
As a result, $(\vc{x}^*, z^*)$ and $w^*$ cannot be returned by the algorithm as the output, a contradiction.
(ii) Assume that $\tilde{\vc{x}}$ is not added as a partial assignment to set $\hat{\mc{X}}$ because the algorithm terminates before such a partial assignment is reached in line 22.
This implies that there must be a partial assignment $\hat{\vc{x}} \in \hat{\mc{X}}$ with $\hat{x}_i = \tilde{x}_i$ for $i=1, \dotsc, j-1$ for some $j \in N$ such that the relaxed DD $\overline{\mc{D}}$ associated with $\mc{M}^C(\hat{\vc{x}})$ for some $C$ is pruned without reaching line 22.
The only possibility for this event is that the length $\bar{w}$ of the longest $r$-$t$ path of $\overline{\mc{D}}$ must satisfy  $\bar{w} \leq w^*$, violating the condition in line 17.
However, because of the facts that $\overline{\mc{D}}$ is a relaxed DD associated with $\mc{M}^C(\hat{\vc{x}})$, and that $(\tilde{\vc{x}}, \tilde{z})$ must be a feasible solution to $\mc{M}^C(\hat{\vc{x}})$, we conclude that $\bar{w} \geq \tilde{w}$.
Combining this inequality with that of the line above, we obtain $\tilde{w} \leq w^*$, a contradiction to the initial assumption on the optimal value of the problem.
(iii) Assume that $\tilde{\vc{x}}$ is not added as a partial assignment to set $\hat{\mc{X}}$ because the longest $r$-$u$ path chosen in line 22 deviates from $\tilde{\vc{x}}$ for some node $u$ of the the path associated with $\tilde{\vc{x}}$.
This implies that there must be a partial assignment $\hat{\vc{x}} \in \hat{\mc{X}}$ with $\hat{x}_i = \tilde{x}_i$ for $i=1, \dotsc, j-1$ for some $j \in N$ such that the relaxed DD $\overline{\mc{D}}$ associated with $\mc{M}^C(\hat{\vc{x}})$ for some $C$ has an exact cut set that contains node $u$ in some layer $k \geq j$.
Let $\dot{\vc{x}}$ be the solution encoding the longest $r$-$u$ path of $\overline{\mc{D}}$ that is chosen in line 22 in place of $\tilde{\vc{x}}$.
It follows from definition of exact cut set that an extension of $\dot{\vc{x}}$ with components $\dot{x}_i = \tilde{x}_i$ for $i = k, \dotsc, n$ and $\dot{z} = \tilde{z}$ is a feasible solution to $\mc{H}$. 
However, the above assumption implies that the length $\dot{w}$ of the path encoded by $(\dot{\vc{x}}, \dot{z})$ is no less than $\tilde{w}$.
If $\dot{w} > \tilde{w}$, then we have reached a contradiction to the assumption that $\tilde{w}$ is the optimal value of $\mc{H}$.
If $\dot{w} = \tilde{w}$, we may repeat the contradiction arguments for the new point $(\dot{\vc{x}}, \dot{z})$ until we exhaust all possible replacements for the assumed optimal solution.
The last such solution must either fall in case (i) or case (ii) above, reaching contradiction.
\Halmos
\endproof
\medskip
\proof{\textbf{Proof of Theorem~\ref{thm:DD master equivalence}}.}
First, we give a useful observation from Algorithm~\ref{alg:DD master}.
It follows from this algorithm that state values $(s^+_u, s^-_u)$ at each node $u \in \mc{U}_j$ records the number of time periods passed since the last start-up and shut-down of the unit at time $j$, as these values reset to 1 whenever there is change in the unit status over two consecutive periods, and are incremented by one if the unit status remains the same.
These state values imply the status of the unit at the beginning of each period, i.e., the unit is down when $s^+_u \geq s^-_u$, and it is up otherwise.

For the direct implication, assume that $(\dot{\vc{x}}, \dot{z})$ encodes an $r$-$t$ path of length $\dot{c}$ in the equivalence class formed by $\mc{D}$. 
We construct an extended solution $(\dot{\vc{x}}, \dot{\vc{y}}, \dot{\bar{\vc{y}}}, \dot{\vc{q}}, \dot{z})$ of \eqref{eq:classic BD master} with objective value $\dot{c}$ as follows.
First, we show that $\dot{x}_j \in \{0, 1\}$ for all $j \in \mc{T}$ and $\dot{z} \in [-\Gamma, \Gamma]$, thereby satisfying domain constraints in \eqref{eq:classic BD master}.
The definition of equivalence class implies that, for each $j \in \mc{T}$, there exists a node pair $(u, v) \in \mc{A}_j \times \mc{A}_{j+1}$ such that $\dot{x}_j \in [l_{(u,v)}^{\text{min}}, l_{(u,v)}^{\text{max}}]$, i.e., the variable value belongs to the interval defined by the min and max label values of the arcs connecting nodes $u$ and $v$; see Section~\ref{subsec:equivalence}.
It follows from the construction of $\mc{D}$ in Algorithm~\ref{alg:DD master} that two consecutive nodes in layers $1, \dotsc, T$ cannot be connected with multiple arcs with different label values, since different arc labels lead to different state values for the head node.
As a result, we must have $l_{(u,v)}^{\text{min}} = l_{(u,v)}^{\text{max}} \in \{0,1\}$.
The argument for $\dot{z} \in [-\Gamma, \Gamma]$ follows directly from the construction and label values assigned to the arcs at the last layer of $\mc{D}$.
To construct the extended solution, for each $j \in \mc{T}$, define $\dot{y}_j = 1$ if $\dot{x}_{j-1} = 0$ and $\dot{x}_j = 1$, and $\dot{y}_j = 0$ otherwise.
Similarly, define $\dot{\bar{y}}_j = 1$ if $\dot{x}_{j-1} = 1$ and $\dot{x}_j = 0$, and $\dot{\bar{y}}_j = 0$ otherwise.
These definitions guarantee the satisfaction of \eqref{const:logical}.
Further, let $u \in \mc{U}_j$ be the node at layer $j$ of the path encoding $(\dot{\vc{x}}, \dot{z})$, and define $\dot{q}_j = K_{s^-_u}$ if $\dot{x}_j = 1$ and $\dot{x}_{j-1} = 0$, and $\dot{q}_j = 0$ otherwise.
These value assignments satisfy \eqref{const:expo_cost} because this constraint implies that when a unit changes status from down to up, i.e., $x_j = 1$ and $x_{j-1} = 0$, then $q_j \geq \max_{k = 1,\dotsc,s^-_u} K_k$, as $s^-_u$ represents the number of periods that the unit has been down consecutively before going up.
Since the start-up cost function is assumed to be \blue{logarithmic}, we have that $\max_{k = 1,\dotsc,s^-_u} K_k = K_{s^-_u}$.
Further, since $\dot{q}_j$ takes the maximum value at equality, it yields a non-redundant solution.
It remains to show that the constructed point satisfies constraints~\eqref{const:min_up} and \eqref{const:min_down}.
Assume by contradiction that there exists $j^* \in \mc{T}$ for which constraint~\eqref{const:min_up} is violated by $(\dot{\vc{x}}, \dot{\vc{y}}, \dot{\bar{\vc{y}}}, \dot{\vc{q}}, \dot{z})$. 
Note that \eqref{const:min_up} models two requirements at time $j^*$: (i) if the unit is down, then it could not have started up in the last $L$ time periods; and (ii) if the unit is up, then it could not have started up more than once in the last $L$ time periods.
For the contradiction, assume first that condition (i) is violated, i.e., $\dot{x}_{j^*} = 0$ and there exists $\bar{j} \in \{j^*-L+1, \dotsc, j^*\}$ such that $\dot{y}_{\bar{j}} = 1$.  
It follows from construction that $\dot{x}_{\bar{j}} = 1$ and $\dot{x}_{\bar{j}-1} = 0$.
Let $(u,v) \in \mc{U}_{\bar{j}} \times \mc{U}_{\bar{j}+1}$ be the nodes at layers $\bar{j}$ and $\bar{j}+1$ of the $r$-$t$ path encoding $(\dot{\vc{x}}, \dot{z})$.
We have $s^+_u \geq s^-_u$ because of the earlier argument on the state values.
Since the arc connecting $u$ to $v$ on the $r$-$t$ path is a 1-arc ($\dot{x}_{\bar{j}} = 1$), we have $s^+_v = 1$; see line 6 of Algorithm~\ref{alg:DD master}.
It follows from the algorithm steps that the only way for the unit to be able to shut down after $\bar{j}$ is to satisfy the condition of line 9, i.e., $s^+_h \geq L$ for some node $h$ in layers $\bar{j}+L, \dotsc, T$.
Since $\bar{j} \geq j^*-L+1$ by definition, we obtain that $j^* < \bar{j}+L$, and hence $\dot{x}_{j^*}$ cannot be equal to zero, a contradiction.
Next, assume that condition (ii) above is violated for the contradiction assumption, i.e., the unit starts up more than once in periods $j^*-L+1, \dotsc, j^*$.
It is easy to verify that there exists a layer $\tilde{j} \in \{j^*-L+1, \dotsc, j^*\}$ for which condition (i) is violated.
Therefore, an argument similar to that of condition (i) yields the desired contradiction.
The contradiction for violating constraint~\eqref{const:min_down} is obtained similarly due to the symmetry in the problem structure.
We conclude that $(\dot{\vc{x}}, \dot{\vc{y}}, \dot{\bar{\vc{y}}}, \dot{\vc{q}}, \dot{z})$ is feasible to \eqref{eq:classic BD master}.
We next show that the length of the $r$-$t$ path encoding $(\dot{\vc{x}}, \dot{z})$ matches the objective value of $(\dot{\vc{x}}, \dot{\vc{y}}, \dot{\bar{\vc{y}}}, \dot{\vc{q}}, \dot{z})$.
The proof follows from considering the contribution of three terms in the objective function of \eqref{eq:classic BD master}.
First, the contribution of each variable assignment $\dot{x}_j = 1$ to the objective function is $c_f$, which is captured in the weight of the associated 1-arcs of the $r$-$t$ path through Algorithm~\ref{alg:DD master} in lines 6 and 8.
Second, the contribution of the start-up status of the unit to the objective function is $\dot{q}_j = K^1_{s^-_u}$ where $u$ is the node at layer $j$ of the path encoding $(\dot{\vc{x}}, \dot{z})$, which is captured in line 6 of Algorithm~\ref{alg:DD master}.
Third, the contribution of $\dot{z}$ in the objective function is directly considered in the arc weight in the last layer of the equivalence class of $\mc{D}$.

For the reverse implication, assume that $(\dot{\vc{x}}, \dot{\vc{y}}, \dot{\bar{\vc{y}}}, \dot{\vc{q}}, \dot{z})$ is a non-redundant solution of \eqref{eq:classic BD master} with objective value $\dot{c}$.
We show that the projected point $(\dot{\vc{x}}, \dot{z})$ encodes an $r$-$t$ path of length $\dot{c}$ in the equivalence class formed by $\mc{D}$. 
Assume by contradiction that no such path exists in the equivalence class of $\mc{D}$.
Therefore, there exists an $r$-$t$ path $\mathtt{P}$ of $\mc{D}$ with associated point $(\bar{\vc{x}},\bar{z})$ and a layer number $j^* \in \mc{T}\setminus\{T\}$ such that $\bar{x}_j = \dot{x}_j$ for $j \in \{1, \dotsc, j^*\}$ and the node $u$ at layer $j^*+1$ of $\mathtt{P}$ does not have an outgoing arc with label value $\dot{x}_{j^*+1}$.
There are four cases.
For the first case, assume that $\dot{x}_{j^*} = 0$ and $\dot{x}_{j^*+1} = 0$.
It follows from the state value definitions that $s^+_u \geq s^-_u$ as $\dot{x}_{j^*} = \bar{x}_{j^*} = 0$. 
Line 4 of Algorithm~\ref{alg:DD master} implies that $u$ has an outgoing arc with label value $\dot{x}_{j^*+1} = 0$, a contradiction.
For the second case, assume that $\dot{x}_{j^*} = 1$ and $\dot{x}_{j^*+1} = 0$.
Because of constraints~\eqref{const:logical} and \eqref{const:min_up}, we must have $\dot{x}_j = \bar{x}_j = 1$ for $j = j^*-L+1, \dotsc, j^*$.
Let $v$ be the node at layer $j^*-L+1$ of $\mathtt{P}$.
Since $\bar{x}_j = 1$ for $j = j^*-L+1, \dotsc, j^*$, it follows from line 8 of Algorithm~\ref{alg:DD master} that $s^+_u = s^+_v + L - 1 \geq L$ as $s^+_v \geq 1$.
As a result, the condition of line 9 of Algorithm~\ref{alg:DD master} is satisfied at node $u$, which leads to creating a 0-arc as an outgoing arc of $u$, a contradiction.
For the third case, where $\dot{x}_{j^*} = 1$ and $\dot{x}_{j^*+1} = 1$, the contradiction is obtained similarly to the first case due to the symmetry of the problem.
For the fourth case, where $\dot{x}_{j^*} = 0$ and $\dot{x}_{j^*+1} = 1$, the contradiction is achieved similarly to the second case due to the symmetry of the problem.
Finally, since $\dot{z} \in [-\Gamma, \Gamma]$, it must belong to the equivalence class of $\mc{D}$ as $-\Gamma$ and $\Gamma$ are used as label values of the pair of arcs connecting every node of layer $T+1$ to the terminal node in $\mc{D}$.
We conclude that $(\dot{\vc{x}}, \dot{z})$ encodes an $r$-$t$ path of length $\dot{c}$ in the equivalence class formed by $\mc{D}$.
To show that the objective value $\dot{c}$ is equal to the length of the $r$-$t$ path encoding $(\dot{\vc{x}}, \dot{z})$, we note that $\dot{q}^1_j = K^1_{s^-_h}$ where $h$ is the node at layer $j$ of the $r$-$t$ path, since otherwise $(\dot{\vc{x}}, \dot{\vc{y}}, \dot{\bar{\vc{y}}}, \dot{\vc{q}}, \dot{z})$ would be a redundant solution of \eqref{eq:classic BD master}.
The rest of the proof follows from an argument similar to that of the direct implication case given above.
\Halmos
\endproof
\medskip
\proof{\textbf{Proof of Proposition~\ref{prop:relaxed DD UCP}}.}
We need to show that for each path $\bar{\mathtt{P}}$ of $\bar{\mc{D}}$ with encoding point $(\bar{\vc{x}},\bar{z})$, there exists a path $\tilde{\mathtt{P}}$ of $\tilde{\mc{D}}$ with encoding point $(\tilde{\vc{x}},\tilde{z})$ such that $\bar{\vc{x}} = \tilde{\vc{x}}$, $\bar{z} = \tilde{z}$, and $w(\bar{\mathtt{P}}) \geq w(\tilde{\mathtt{P}})$, where $w(.)$ represent the length of the path.
Define $\bar{\mt{P}}^k$ to be the path composed of the first $k$ arcs of $\bar{\mt{P}}$ for $k \in \mc{T}$.
We prove the result by induction on the size of $k$, i.e., we show that for each $\bar{\mathtt{P}}^k$ with encoding point $\bar{\vc{x}}^k$, there exists a path $\tilde{\mathtt{P}}^k$ of $\tilde{\mc{D}}$ with encoding point $\tilde{\vc{x}}^k$ such that $\bar{\vc{x}}^k = \tilde{\vc{x}}^k$ and $w(\bar{\mathtt{P}}^k) \geq w(\tilde{\mathtt{P}}^k)$.
Note that it suffices to show the above result for $k$ up to $T$, since each node at node layer $T+1$ is connected to the terminal node by two arcs with labels $-\Gamma$ and $\Gamma$ in both $\bar{\mc{D}}$ and $\tilde{\mc{D}}$, and therefore a matching arc can always be found for the desired path.
For the base of induction, i.e., $k=1$, it follows from Algorithm~\ref{alg:DD master} that $\bar{\mc{D}}$ has two arcs with labels $0$ and $1$ going out of the root node.
These two arcs remain in $\tilde{\mc{D}}$ with the same weights according to Definition~\ref{def:modified alg}.
For the induction hypothesis, assume that, for some $k \geq 2$, there exists a path $\tilde{\mathtt{P}}^{k-1}$ of $\tilde{\mc{D}}$ with encoding point $\tilde{\vc{x}}^{k-1}$ such that $\bar{\vc{x}}^{k-1} = \tilde{\vc{x}}^{k-1}$ and $w(\bar{\mathtt{P}}^{k-1}) \geq w(\tilde{\mathtt{P}}^{k-1})$.
We show that there exists a path $\tilde{\mathtt{P}}^k$ of $\tilde{\mc{D}}$ with encoding point $\tilde{\vc{x}}^k$ such that $\bar{\vc{x}}^k = \tilde{\vc{x}}^k$ and $w(\bar{\mathtt{P}}^k) \geq w(\tilde{\mathtt{P}}^k)$.
There are four cases.
(i) Assume that $\bar{x}^k_{k-1} = 0$ and $\bar{x}^k_{k} = 0$.
Note that $\bar{x}^k_{k-1} = \bar{x}^{k-1}_{k-1} = 0$ as they are defined over the same path $\bar{\mathtt{P}}$.
It follows from the induction hypothesis that $\bar{x}^{k-1}_{k-1} = \tilde{x}^{k-1}_{k-1} = 0$.
Let $u$ be the node at node layer $k$ of the path $\tilde{\mathtt{P}}^{k-1}$.
We consider two cases for this node.
For the first case, assume that $u$ is not a merged node.
As a result, it has been created directly from the modified algorithm of Definition~\ref{def:modified alg}, which implies that $s^+_u \geq s^-_u$.
Following line 4 of Algorithm~\ref{alg:DD master}, there is a 0-arc going out of $u$ in $\tilde{\mc{D}}$.
Adding this arc to $\tilde{\mathtt{P}}^{k-1}$ yields the desired path $\tilde{\mathtt{P}}^{k}$, as the contribution to the objective function is 0 for this variable assignment in both paths $\bar{\mathtt{P}}^{k}$ and $\tilde{\mathtt{P}}^{k}$.
For the second case, assume that $u$ is a merged node.
Let $\{v_1, \dotsc, v_r\}$ be the set of nodes that have been merged into $u$ using the merging operation of Definition~\ref{def:merging}.
Since $u$ has an incoming arc with label 0, it means that one of the nodes $v_i$ in the above set must have an incoming arc with label 0.
The modified algorithm implies that $s^+_{v_i} \geq s^-_{v_i}$.
The merging condition, then, dictates that $s^+_{v_i} \geq s^-_{v_i}$ for all $i \in \{1, \dotsc, r\}$.
It follows from the merging equations that $s^+_u = \max_i s^+_{v_i} \geq \max_i s^-_{v_i} = s^-_u$.
Therefore, an argument similar to that of the first case yields the desired result.
(ii) Assume that $\bar{x}^k_{k-1} = 0$ and $\bar{x}^k_{k} = 1$.
Let $h$ be the node at node layer $k$ of $\bar{\mathtt{P}}^{k}$.
It follows from line 5 of Algorithm~\ref{alg:DD master} that $s^-_h \geq \ell$ and that the arc weight at layer $k$ is $c_f + K_{s^-_h}$.
It is clear from the definition of $s^-_h$ that $\bar{x}^k_{k-s^-_h-1} = 1$ and $\bar{x}^k_{k-s^-_h -1 + i} = 0$ for all $i = 1, \dotsc, s^-_h$.
It follows from the induction hypothesis that $\tilde{x}^{k-1}_{k-s^-_h-1} = 1$ and $\tilde{x}^{k-1}_{k-s^-_h -1 + i} = 0$ for all $i = 1, \dotsc, s^-_h$.
Let $u$ and $v$ be the nodes at node layers $k-s^-_h+1$ and $k$ of $\tilde{\mathtt{P}}^{k-1}$, respectively.
We consider two cases for node $u$.
For the first case, assume that $u$ is not a merged node.
As a result, it has been created directly from the modified algorithm of Definition~\ref{def:modified alg}, which implies that $s^-_u = s^=_u = 1$.
For the second case, assume that $u$ is a merged node.
Let $\{u_1, \dotsc, u_r\}$ be the set of nodes that have been merged into $u$ using the merging operation of Definition~\ref{def:merging}.
Since $u$ has an incoming arc with label 0 right after an arc with label 1, it implies that one of the nodes $u_i$ in the above set must have $s^-_{u_i} = s^=_{u_i} = 1$.
It follows from the modified algorithm that $s^-_{u} \geq 1$ and $s^=_u \leq 1$, which holds for the first case as well.
Using an argument similar to that given above for the nodes at layers $k-s^-_h+1$ to $k$ of $\tilde{\mathtt{P}}^{k-1}$, we conclude that $s^-_v \geq s^-_u + s^-_h - 1 \geq s^-_h $ and $s^=_v \leq s^=_u + s^-_h - 1 \leq s^=_h$, where the first inequalities follow from the state definitions through merging operation applied on successive node layers, and the second inequalities hold because of the relations above on state values at node $u$.
As a result, the condition in line 5 of Algorithm~\ref{alg:DD master} is satisfied at node $v$ as $s^-_v \geq s^-_h \geq \ell$, meaning that it has an outgoing arc with label 1 and weight $c_f + K_{s^=_v} \leq c_f + K_{s^=_h}$.
Adding this arc to $\tilde{\mathtt{P}}^{k-1}$ yields the desired path $\tilde{\mathtt{P}}^{k}$.
(iii) Assume that $\bar{x}^k_{k-1} = 1$ and $\bar{x}^k_{k} = 1$.
The result is obtained by using an argument similar to that of case (i) due to symmetry.
(iv) Assume that $\bar{x}^k_{k-1} = 1$ and $\bar{x}^k_{k} = 0$.
An argument similar to that of case (ii) proves the result due to symmetry. \Halmos
\endproof
\medskip
\proof{\textbf{Proof of Proposition~\ref{prop:ramp_equivalent}}.}
We prove the equivalence by showing that for each feasible variable assignment, both sets of constraints impose the same restrictions.
There are four cases for each $i \in N$ and $j \in \mc{T}$.
(i) Assume that $x^i_{j-1}=x^i_j=0$. 
It follows from constraints \eqref{const:logical} and \eqref{const:min_up} that $y^i_j=\bar{y}^i_j =0$. 
In this case, the right-hand-side (RHS) of both inequalities \eqref{const:ramp_up} and \eqref{eq:ramp_up_2}, as well as \eqref{const:ramp_down} and \eqref{eq:ramp_down_2} reduce to zero. (ii) Assume that $x^i_{j-1}=0$ and $x^i_j=1$. 
It follows from constraint \eqref{const:logical} that $y^i_j=1$ and $\bar{y}^i_j = 0$. 
In this case, the RHS of both constraints \eqref{const:ramp_up} and \eqref{eq:ramp_up_2} reduce to $\SU^i$.
Further, inequality \eqref{const:pow_gen} implies that $p^i_{j-1} = 0$ as $x^i_{j-1}=0$.
As a result, both constraints \eqref{const:ramp_down} and \eqref{eq:ramp_down_2} become redundant in the description of $\mc{E}$ and $\mc{G}$, since their LHS is a non-positive quantity $-p^i_j \leq 0$ and their RHS are non-negative values as $\RD^i \geq 0$ and $\RD^i - \SD^i \geq 0$ by assumption.
(iii) Assume that $x^i_{j-1}=1$ and $x^i_j=0$. 
It follows from constraint \eqref{const:logical} that $y^i_j=0$ and $\bar{y}^i_j = 1$.
In this case, the RHS of both constraints \eqref{const:ramp_down} and \eqref{eq:ramp_down_2} reduce to $\SD^i$.
Further, inequality \eqref{const:pow_gen} implies that $p^i_{j} = 0$ as $x^i_{j}=0$.
As a result, both constraints \eqref{const:ramp_up} and \eqref{eq:ramp_up_2} become redundant in the description of $\mc{E}$ and $\mc{G}$, since their LHS is a non-positive quantity $-p^i_{j-1} \leq 0$ and their RHS are non-negative values as $\RU^i \geq 0$ and $\RU^i - \SU^i \geq 0$ by assumption.
(iv) Assume that $x^i_{j-1}=x^i_j=1$. 
It follows from constraints \eqref{const:logical} and \eqref{const:min_down} that $y^i_j=\bar{y}^i_j =0$. 
In this case, the RHS of both \eqref{const:ramp_up} and \eqref{eq:ramp_up_2} reduce to $\RU^i$, and the RHS of both \eqref{const:ramp_down} and \eqref{eq:ramp_down_2} reduce to $\RD^i$. \Halmos     
\endproof

\section{Comparison of DD-BD Method with the Literature} \label{sec:bd_literature}

In the literature, DDs have been used in conjunction with BD framework in different contexts.
In this section, we give a detailed comparison between DD-BD approach of Section~\ref{sec:Benders} and those used in the literature.
In particular, we show that our framework generalizes other existing approaches, while addressing their shortcomings.


\cite{van2017decision} proposes a BD framework for linear mixed integer programs that uses a similar approach to ours in defining the master problem over integer variables, and solving the subproblems for continuous variables.
A DD is constructed to represent the solution set of the master problem, and is successively refined through cuts obtained from the subproblems.
There is, however, a fundamental difference between our approach and that of~\cite{van2017decision} as they use separation for the feasibility cuts only since their proposed DD cannot contain any continuous variable.
To incorporate the optimality cuts, the authors propose a so-called ``cost-tuple" approach where, for each optimality cut, a cost parameter is calculated at each node of the DD to record the contribution of the variable assignments in relation to that cut.
Such an approach requires the underlying DD to be \textit{non-reduced}, limiting its practicality.
In contrast, our proposed DD-BD approach provides a unified framework that treats both feasibility and optimality cuts similarly through separation as the continuous variable is directly incorporated in the DD layers.
To elaborate, we first give a brief review of the cost-tuple approach used in~\cite{van2017decision} to handle optimality cuts. 

The cost-tuple approach is given in Algorithm~\ref{alg:cost-tuple} where, without loss of generality, the master problem of a MIP is defined as $\max \{ z \, | \, \vc{x} \in \mc{P} \}$ with a bounded $\mc{P} \subseteq \Z^n$. 
Consider $\mc{D} = (\mc{U},\mc{A},l(\cdot))$ to be the DD representing $\mc{P}$, with a property that each node (except the terminal) has a unique incoming arc, and the terminal receives a unique arc from each node of the previous layer.
Let inequalities of the form $z \leq \vc{\alpha}^j \vc{x} + \alpha_0^j$, for $j \in J$, represent the set of added optimality cuts at the current iteration of the BD algorithm.
For each node $u \in \mc{U}$, parameter $r^j(u)$ records the reward (cost, for a minimization variant) contribution of variable assignments at that node to the right-hand-side value of cut $z \leq \vc{\alpha}^j \vc{x} + \alpha_0^j$ for $j \in J$.
In this algorithm, the accumulated reward $r(u)$ captures the contribution amount at the intersection of optimality cuts for each unique value assignment of variables, while the final reward $r^*$ computes the optimal value of $\max \{ z \, | \, \vc{x} \in \mc{P}, z \leq \vc{\alpha}^j \vc{x} + \alpha_0^j, \forall j \in J \}$.

\begin{algorithm}[!htbp]                    				  	
	\caption{The cost-tuple approach of~\cite{van2017decision} to handle optimality cuts}          
	\label{alg:cost-tuple}			                        
		\KwData{a DD $\mc{D} = (\mc{U},\mc{A},l(.))$ and a set of optimality cuts $z \leq \vc{\alpha}^j \vc{x} + \alpha_0^j$, for $j \in J$} 
		initialize the root node reward as $r^j(s) = \alpha_0^j$ for $j \in J$
		
		\ForAll{cuts $j \in J$, node layers $i \in N\setminus\{n\}$, and node pairs $(u,v) \in \mc{U}_i \times \mc{U}_{i+1}$}{
		    compute $r^j(v) = r^j(u) + \alpha^j_i l_{a}$ where $a \in \mc{A}_i$ is the unique arc connecting $u$ to $v$
		}
		\ForAll{nodes $u \in \mc{U}_n$}{
		    compute accumulated reward $r(u) = \min_{j \in J}\{r^j(u) + \alpha^j_n l_a\}$ where $a \in \mc{A}_n$ is the unique arc connecting $u$ to the terminal $t$
		}
		compute the final reward $r^* = \max_{u \in \mc{U}_n}r_u$
\end{algorithm}

The main shortcoming of Algorithm~\ref{alg:cost-tuple} is due to the assumption that the DD nodes have a unique incoming arcs.
This requirement leads to an exponential growth of the DD size to distinguish between each $r$-$t$ path. 
As a consequence, it is mentioned in~\cite{van2017decision} that the cost-tuple algorithm is not applicable to reduced DDs.
Reduction in DDs is referred to the operation of merging nodes that share a similar subtree; hence it is regarded as a crucial component in building efficient DDs for practical applications and keeping them within a width limit.
The DD-BD approach proposed in the present paper, on the other hand, can be applied to reduced DDs, and the separation of merged nodes, if needed, will be performed naturally through the separation operation. 
We illustrate this difference in the following example.

\begin{example} \label{ex:BD}
	Consider the final iteration of a BD algorithm with a master problem $\max \{ z \, | \, \vc{x} \in \mc{P} \}$ where $\mc{P} \subseteq \Z^n$ is represented by DD $\mc{D}_1$ of Figure~\ref{fig:bd1}.
	Assume that the following two optimality cuts are added 
	\begin{subequations}
		\begin{align}
	&z \leq 3x_1 + 2x_2 \label{eq:c1} \\
	&z \leq -3x_1 - 5x_2 + 3 \label{eq:c2}
	\end{align}
		\end{subequations}
	It is clear that both points $(0, 0)$ and $(1, 0)$ give an optimal solution with optimal value $z^* = 0$.
	Since $\mc{D}_1$ is a reduced DD, the approach of \cite{van2017decision} cannot be directly applied to obtain the optimal solution.
	It is suggested in that work that, when merging occurs, the component-wise maximum of the reward tuple for merging nodes must be selected as the reward tuple of the merged node.
	This operation guarantees that the final reward value is an upper-bound for the optimal reward value.
	Applying this operation here, however, creates an infinite loop without converging to a solution, as demonstrated next.
	We record the reward value for inequalities \eqref{eq:c1} and \eqref{eq:c2}, respectively, in a tuple $(r^1(u), r^2(u))$ for each node $u$.
	For the root node, we have the reward tuple $(0, 3)$ according to Algorithm~\ref{alg:cost-tuple}.
	For the path including the 1-arc at the first layer, the reward of the end-node is $(3, 0)$.
	Similarly, the path including the 0-arc at the first layer has the end-node reward $(0, 3)$.
	The component-wise maximum yields the reward of $(3, 3)$ for the merged node in the middle.
	The reward at the terminal node is $(3, 3)$, giving the optimal value $\bar{z} = 3$.
	While $\bar{z}$ is an upper-bound for $z^*$, it can already be separated by the current optimality cuts \eqref{eq:c1} and \eqref{eq:c2} for both feasible solutions.
	As a result, the BD approach cannot refine it any further and is stuck at this iteration.
	
	For the DD-BD approach of Algorithm~\ref{alg:DD-BD}, we represent the master problem as $\mc{M} = \max \{ z \, | \, \vc{x} \in \mc{P}, z \in [-M, M] \}$ where valid bounds are imposed on $z$ variable.
	A restricted DD $\mc{D}_2$ can be created for $\mc{M}$ using the equivalence class argument of Section~\ref{subsec:equivalence} as given in Figure~\ref{fig:bd2}.
	Next, we may simply refine the DD through a sequential separation with respect to optimality cuts \eqref{eq:c1} and \eqref{eq:c2}, as shown in Figures~\ref{fig:co1} and \ref{fig:co2}.
	It is easy to verify that a longest path in the resulting DD gives either point $(0, 0)$ or $(1, 0)$ with the optimal value $z^* = 0$.
	Since these solutions give the optimal value of the problem, Algorithm~\ref{alg:DD-BD} is terminated due to bound pruning conditions in line 17. 
\end{example}

\begin{figure}[!hbt]
	\centering
	\begin{subfigure}[b]{0.45\linewidth} 
		\centering
		\includegraphics[scale=0.4]{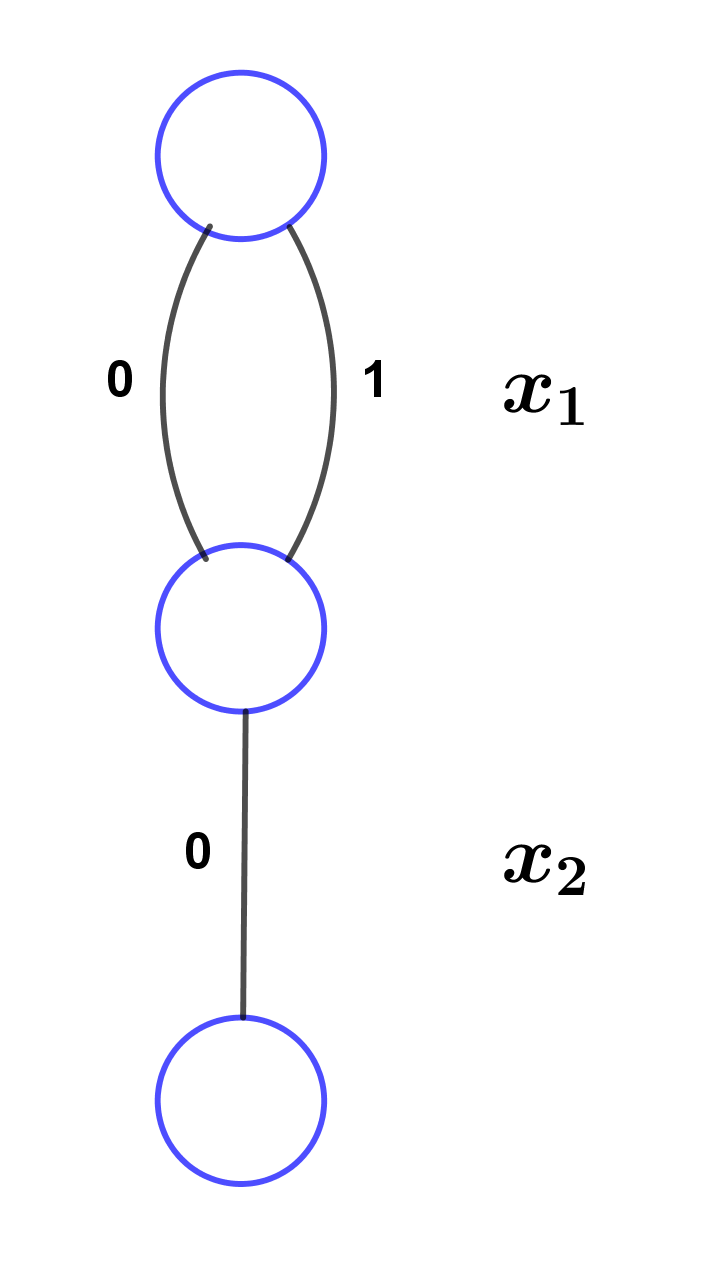}
		\caption{$\mc{D}_1$ without continuous variable}
		\label{fig:bd1}
	\end{subfigure}
	\hspace{0.05\linewidth}
	\begin{subfigure}[b]{0.45\linewidth} 
		\centering
		\includegraphics[scale=0.4]{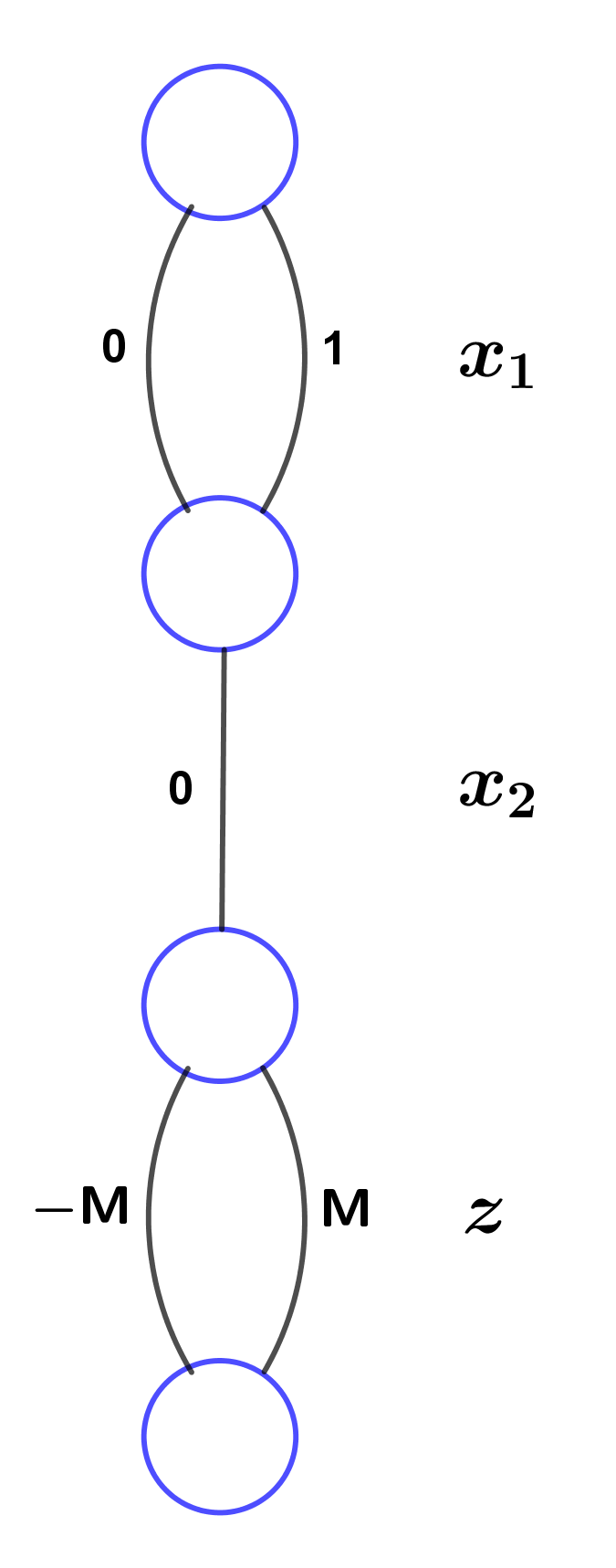}
		\caption{$\mc{D}_2$ with continuous variable}
		\label{fig:bd2}
	\end{subfigure}	
	\caption{DD representation of the master problem of Example~\ref{ex:BD}}
	\label{fig:bd}
\end{figure}

\begin{figure}[!hbt]
	\centering
	\begin{subfigure}[b]{0.45\linewidth} 
		\centering
		\includegraphics[scale=0.4]{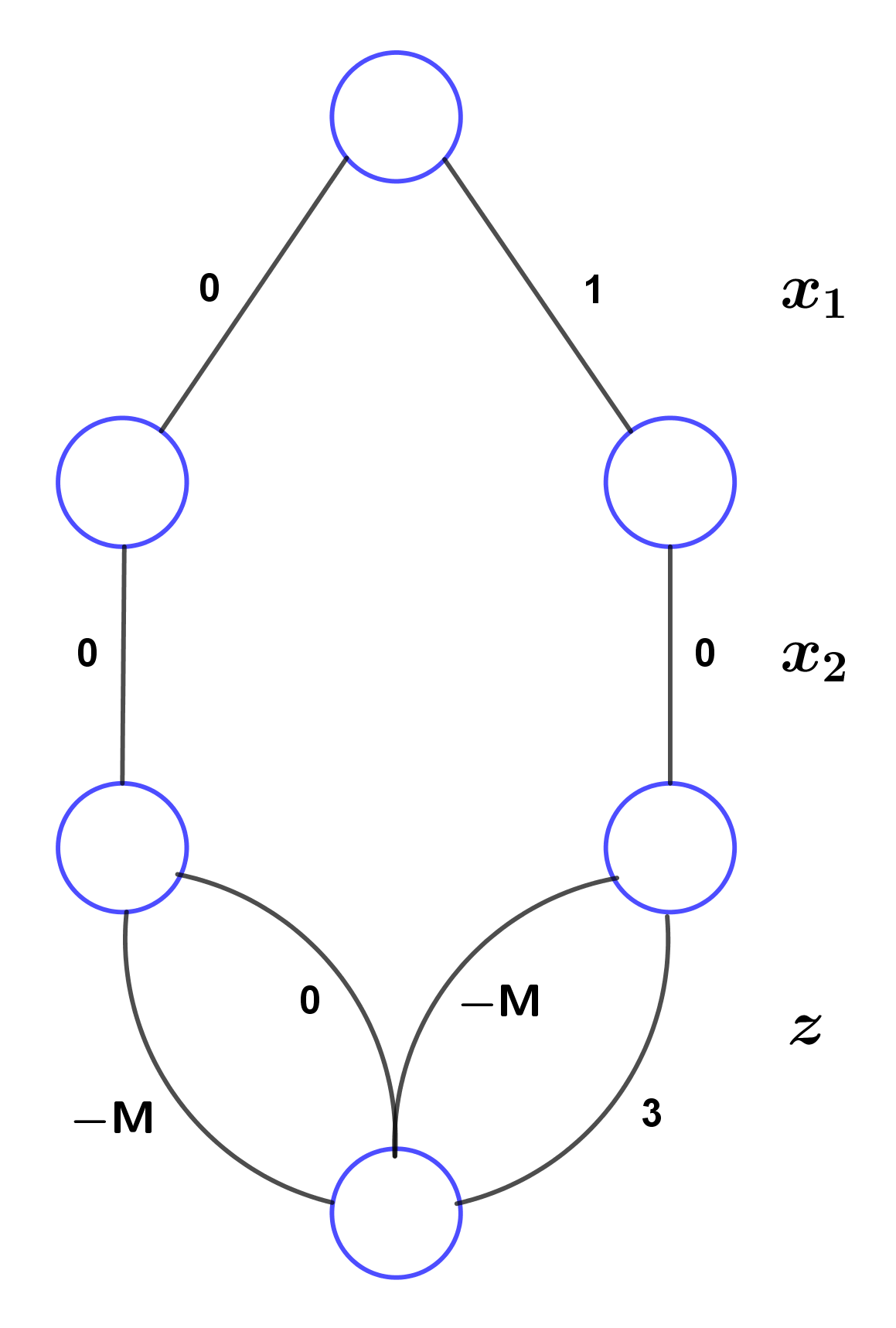}
		\caption{Separation wrt \eqref{eq:c1}}
		\label{fig:co1}
	\end{subfigure}
	\hspace{0.05\linewidth}
	\begin{subfigure}[b]{0.45\linewidth} 
		\centering
		\includegraphics[scale=0.4]{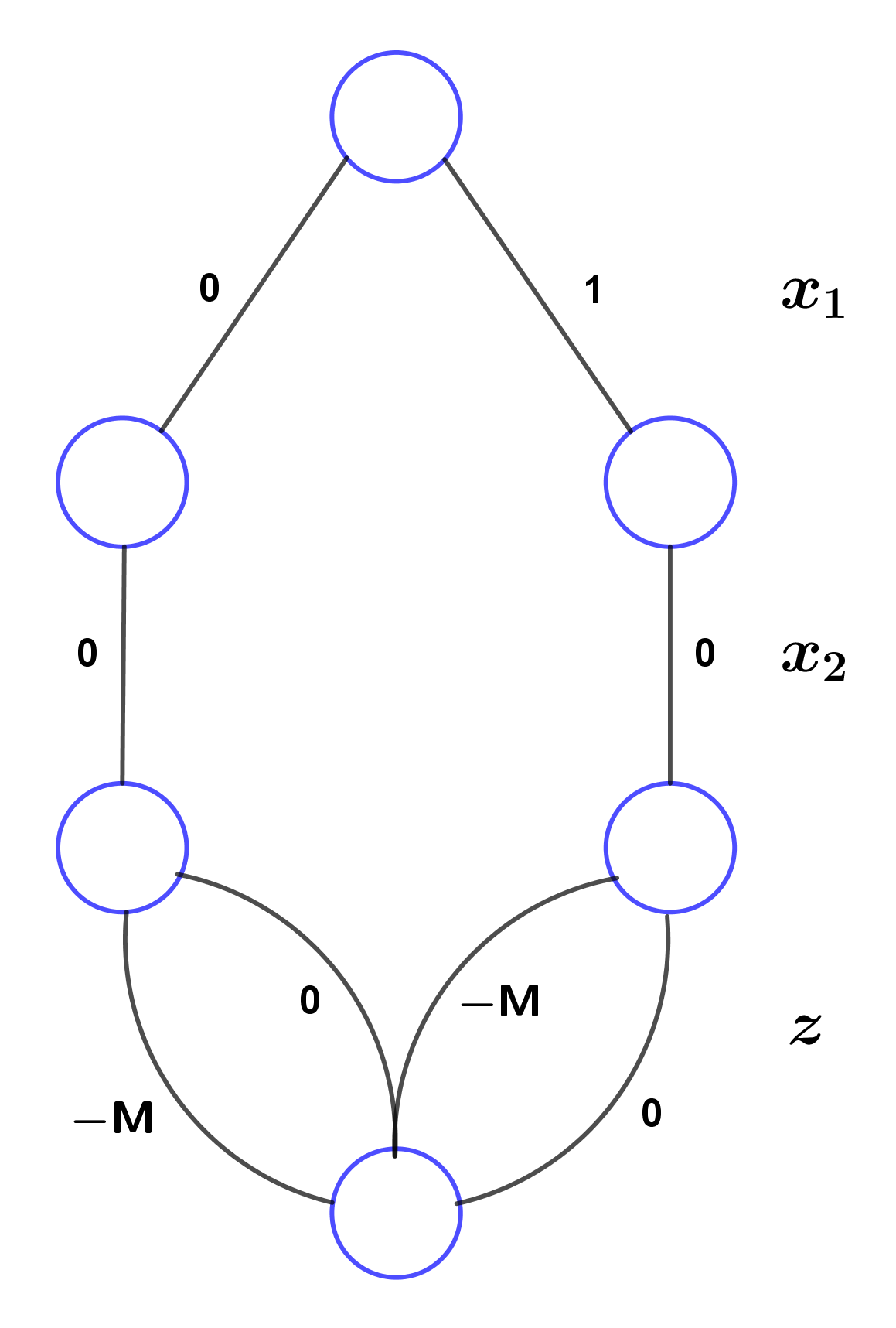}
		\caption{Separation wrt \eqref{eq:c2}}
		\label{fig:co2}
	\end{subfigure}	
	\caption{Separation of optimality cuts for the master problem of Example~\ref{ex:BD}}
	\label{fig:co}
\end{figure}

We conclude this section by showing that the DD-BD approach of Algorithm~\ref{alg:DD-BD} is a generalization of the cost-tuple approach of Algorithm~\ref{alg:cost-tuple}, since the latter can be viewed as a special case of the former when the underlying DD is not reduced.

\begin{proposition} \label{prop:cost-tuple generalization}
Consider the master problem $\max \{ z \, | \, \vc{x} \in \mc{P} \}$ of a MIP with a bounded $\mc{P} \subseteq \Z^n$, and let $\mc{D} = (\mc{U},\mc{A},l(\cdot))$ be the DD representing $\mc{P}$, with a property that each node (except the terminal) has a unique incoming arc, and the terminal receives a unique arc from each node of the previous layer.
Let inequalities of the form $z \leq \vc{\alpha}^j \vc{x} + \alpha_0^j$, for $j \in J$, represent the set of optimality cuts to be added at the current iteration of the BD algorithm.
Define $\bar{\mc{D}} = (\bar{\mc{U}},\bar{\mc{A}},\bar{l}(\cdot))$ to be the DD constructed from $\mc{D}$ by adding a layer representing $z$ variable as the last arc layer.
Then, the reward value obtained from Algorithm~\ref{alg:cost-tuple} applied to $\mc{D}$ is equal to the optimal value obtained from Algorithm~\ref{alg:DD-BD} when applied to $\bar{\mc{D}}$.
\end{proposition}

\proof{Proof.} 
We first describe the relation between $\mc{D}$ and $\bar{\mc{D}}$. 
All nodes and arcs from layer one to layer $n-1$ are similar in both DDs, i.e., $\mc{A}_i = \bar{\mc{A}}_i$ for $i \in N \setminus \{n\}$ and $\mc{U}_i = \bar{\mc{U}}_i$ for $i \in N$.
The difference is that $\bar{\mc{D}}$ has an extra layer that represents $z$ variable.
This layer is created as follows.
For each node $u \in \bar{\mc{U}}_n$, we create a node $v \in \bar{\mc{U}}_{n+1}$, and connect them with an arc $a \in \bar{\mc{A}}_n$ with the same label as the unique arc connecting the replica of node $u$ in $\mc{U}_n$ to the terminal node of $\mc{D}$. 
Then, for each node $v \in \bar{\mc{U}}_{n+1}$, we create two arcs with labels $-M$ and $M$ in layer $\bar{\mc{A}}_{n+1}$ that connect $v$ to the terminal node in $\bar{\mc{D}}$.
These label values represent some valid lower and upper bounds for variable $z$.
According to Algorithm~\ref{alg:DD-BD}, $\bar{\mc{D}}$ is refined with respect to the optimality cuts, and then the longest path is found over the refined DD weighted by the objective function rates.
To perform the refinement, the solution set satisfying each optimality cut is modeled by a DD called \textit{threshold diagram}; see \cite{bergman:ci:va:ho:2016} for an exposure to threshold diagrams. 
The last layer of such DDs corresponds to variable $z$ and contains two arcs with labels representing a lower and upper bound for $z$.
Since for each node $v \in \bar{\mc{U}}_{n+1}$, there exists a unique arc-specified path $\mr{P} = (a_1,\dotsc,a_n) \in \bar{\mc{A}}_1 \times \dotsc \times \bar{\mc{A}}_n$ from the root node to $v$, it follows from the assumption that the threshold diagrams have the same structure as $\bar{\mc{D}}$, with the difference that the upper bound label on arcs connecting $v$ to the terminal is computed as $M^j(v) = \sum_{i=1}^n\alpha^j_i l_{a_i} + \alpha_0^j$ for each optimality cut $j \in J$.
As a result, refining $\bar{\mc{D}}$ with respect to the optimality cuts reduces to intersecting label value intervals associated with arcs at the last layer, which are of the form $[-M, M^j(v)]$ for each node $v \in \bar{\mc{U}}_{n+1}$. 
This intersection yields the lower bound label $-M$, and the upper bound label $M(v) = \min_{j \in J}M^j(v)$ for pair of arcs connected to $v$ in the last layer of the refined DD.
Note that $M(v)$ is equal to $r(v)$ as computed in Algorithm~\ref{alg:cost-tuple}.
Finally, to find the longest path on the refined $\bar{\mc{D}}$, all arcs in layers $1$--$n$ are assigned zero weight and arcs in the last layer are assigned weight 1 deduced from the objective function of the master problem.
Therefore, the length of the longest $r$-$t$ path in $\bar{\mc{D}}$ is computed as $w^* = \max_{v \in \bar{\mc{U}}_{n+1}}M(v)$, which is equal to $r^*$ obtained from Algorithm~\ref{alg:cost-tuple}.
Since $\bar{\mc{D}}$ represents an exact DD for the master problem and has been refined with respect to optimality cuts, $w^*$ is returned as the output of Algorithm~\ref{alg:DD-BD}.
\Halmos
\endproof
\end{APPENDICES}

\clearpage
\end{document}